\documentclass[]{aiaa-tc}

 \title{Temporal Stabilisation of Flux Reconstruction on Linear Problems}

 \author{
    Will Trojak%
    \thanks{PhD Candidate, AIAA Student Member, wt247@cam.ac.uk},
  \ Rob Watson%
  	\thanks{Research Fellow, AIAA Member}
  \ and Paul G. Tucker%
  	\thanks{Professor, AIAA Assoc. Fellow}  
  \thanksibid{0}\\
  {\normalsize\itshape
   Department of Engineering, University of Cambridge, Cambridge, UK, CB2 1PZ}\\
 }
 
 \AIAApapernumber{YEAR-NUMBER}
 \AIAAconference{Conference Name, Date, and Location}
 \AIAAcopyright{\AIAAcopyrightD{YEAR}}

 \newcommand{\etal}{\emph{et~al.}}

\newcommand{\hfd}{\hat{f}^{\delta}}
\newcommand{\hf}{\hat{f}}
\newcommand{\hfI}{\hat{f}^{I}}
\newcommand{\hfC}{\hat{f}^{C}}
\newcommand{\hud}{\hat{u}^{\delta}}
\newcommand{\hu}{\hat{u}}

\newcommand\px[2]{\frac{\partial #1}{\partial {#2}}}

\newcommand\dx[2]{\frac{d #1}{d #2}}

\graphicspath{ {Figs/} }

\usepackage{graphicx}
\usepackage{graphics}
\usepackage{upgreek}
\usepackage{float}
\usepackage{epstopdf}
\usepackage{amsmath}
\usepackage{wrapfig}
\usepackage{amssymb}
\usepackage{pdflscape}
\usepackage{listings}
\usepackage{cite}
\usepackage{subcaption}
\usepackage{multirow}
\usepackage[percent]{overpic}
\usepackage{multicol}
\usepackage{color}

\begin{document}

\maketitle

\begin{abstract}
Filtering is often used in Large Eddy Simulation with a global filter width, instead here a filter width in the reference domain of high order Flux Reconstruction is considered. It is shown via Von Neumann analysis how filtering effects the dispersion and dissipation of the scheme when spatially and temporally discretised. With it being shown that filtering stabilises the scheme temporally by upto $25\%$ for forth order FR. The impact of filtering on error production is calculated, highlighting the reduction in convective velocity caused and showing numerically the impact on order of accuracy. Finally, the turbulent Taylor-Green case is used to understand the effect of reference domain filtering on the transition to turbulence, and a filter Reynolds number is defined that is shown to be useful in understanding the effect of filtering on simulations.

\end{abstract}

\section*{Nomenclature}
\begin{multicols}{2}
\begin{tabbing}
  XXXX \= \kill
    \textit{Roman}\\
	$a$ \> convective velocity \\	  
	$c(k)$ \> wavespeed at wavenumber $k$ \\
	$\mathbf{C}_{+1}$ \> downwind cell FR matrix \\
	$\mathbf{C}_0$ \> centre cell FR matrix \\
	$\mathbf{C}_{-1}$ \> upwind cell FR matrix \\
	$\mathbf{D}$ \> first derivative matrix \\
	$\mathbf{e}_n$ \> $n^{\mathrm{th}}$ element solution point error\\ 
	$f$ \> flux variable \\
	$h_l \:\mathrm{\&}\: h_r$ \> left and right correction functions\\
	$J_n$ \> $n^{\mathrm{th}}$ cell Jacobian\\
	$k_{nq}$ \> solution point Nyquist wavenumber, $(p+1)/\delta_j$\\
	$\hat{k}$ \> $k_{nq}$ normalised wavenumber, $[0,\pi]$ \\
	$M_a$ \> Mach number \\ 
	$p$ \> solution polynomial order \\
	$P_r$ \> Prandtl number \\
	$\mathbf{Q}$ \> spatial scheme matrix \\
	$R_e$ \> Reynolds number \\
	$\mathbf{R}$ \> update matrix \\
	$\mathbf{S}$ \> filter matrix \\
	$u$ \> primitive variable \\
	$\mathbf{W}$ \> eigenvector matrix \\	
	
	\textit{Greek}\\
	$\beta_n$ \> weight of $n^{\mathrm{th}}$ mode \\
	$\pmb{\beta}$ \> array of mode weights \\
	$\gamma$ \> grid geometric expansion factor \\
	$\mathbf{\Gamma}$ \> diagonal eigenvalue matrix\\
	$\delta_j$ \> mesh spacing\\	
	$\iota$ \> VCJH correction function parameter \\
	$\mathbf{\Lambda}$ \> normalised diagonal eigenvalue matrix \\
	$\xi$ \> local computation spatial variable \\
	$\sigma$ \> Gaussian filter width \\
	$\varsigma$ \> Gaussian bump solution width \\ 
	$\tau$ \> time step \\
	$\mathbf{\Omega}$ \> solution domain \\
	$\mathbf{\Omega}_n$ \> $n^{\mathrm{th}}$ solution sub-domain \\
	$\hat{\mathbf{\Omega}}$ \> reference sub-domain\\
	
	\textit{Subscript}\\
	$\mathrm{\bullet}_l$ \> variable at left of cell\\
	$\mathrm{\bullet}_r$ \> variable at right of cell\\
	\textit{Superscript}\\
	$\hat{\mathrm{\bullet}}$ \> transformed value \\
	$\overline{\mathrm{\bullet}}$ \> locally fitted polynomial of value \\ 
	$\mathrm{\bullet}^c$ \> common value at interface\\
	$\mathrm{\bullet}^{\delta}$ \> discontinuous value\\
\end{tabbing}
\end{multicols}

\section{Introduction}\label{sec:intro}
	Flux Reconstruction, a high order method originally proposed by Huynh~\cite{Huynh2007}, has proved to be largely robust~\cite{Trojak2017} and is gaining interest as a commercial CFD tool --- particularly as FR seems to overcome some of the computational limitations confronted by lower order methods when performing Large Eddy Simulation~\cite{Vincent2017}. LES itself is of industrial interest as it offers the possiblity of performing without the turbulence modelling uncertainty associated with Reynolds-Averaged approaches, giving improved aerodynamic and aerothermal predictions as a result.

        However, a hurdle that must be confronted by all high order methods is the limitations imposed by the stability of temporal integration, and the effect that this can have on the temporal integration. This can mean that explicit time steps, in the case of FR, are between half and one fifth of those expected from a low order method such as finite volume. These lower CFL limits introduced by FR have the potential to become limiting when running the larger cases that FR makes possible and LES/DNS demand. It has long been known that filtering can aid the temporal stability of schemes, for example Fischer and Mullen~\cite{Fischer2001} showed the stabilising effect of sharp spectral filtering on spectral element methods.	
	
	The authors, in this paper, aim to present the generalised fully discrete von Neumann analysis of Flux Reconstruction and present how, on uniform grids, the harmonic modes of the system impact the stability of the function projection. The use of filtering as a method of instability amelioration is then introduced and the impact on the fully discretised scheme is investigated, along with the effect that reference domain filtering has on the implicit filter of the scheme. A key tool used in the assessment of filtering on both semi- and fully-discrete schemes is an analytical error study, with an extension made here so that the fully discretised version may be considered. Finally, numerical tests are performed to study the effect of filtering when implemented, with the aim of both validating the analytical results and showing how these methods impact the calculation of fully turbulent flows.

\section{Flux Reconstruction}
	As discussed in Section~\ref{sec:intro}, the numerical method under investigation is Flux Reconstruction. We will now give a brief overview of the mechanics of the scheme by presenting a 1D first order example. However, further information on the derivation, stability, and implementation can be found in \cite{Vincent2010,Vincent2015,Trojak2018,Witherden2016a,Jameson2012,Witherden2014,Williams2014a,Williams2014,Jameson2011}.

Let the domain of the solution be defined as $\mathbf{\Omega}$. We then perform domain decomposition such that we have sub-domains, $\mathbf{\Omega}_i$, that fulfil:
	\begin{equation}
		\mathbf{\Omega} = \bigcup_{i=1}^N \mathbf{\Omega}_i \quad \mathrm{and} \quad \bigcap_{i=1}^N \mathbf{\Omega}_i= \emptyset
	\end{equation}

	We then define a reference domain, $\hat{\mathbf{\Omega}} \in [-1,1]$, and associate a Jacobian transformation from one domain to the other such that $\mathbf{J}_i : \mathbf{\Omega}_i \mapsto \hat{\mathbf{\Omega}}_i$. Within the sub-domains and reference sub-domain we position solution points such that $\mathbf{x}_{j,p} = \{x_0,\dots,x_p: x_i \in \mathbf{\Omega}_j\}$ and $\pmb{\xi}_{p} = \{\xi_0,\dots,\xi_p: \xi_i \in \hat{\mathbf{\Omega}}\}$. Collocated points on the edge of the sub-domain are also positioned and in 1D this gives $\pmb{\xi}_f = \{-1,1\}$. The positioning of points is of great importance and was investigated for tetrahedra by Witherden~\cite{Witherden2016}. With the domain of the solution established we now introduce both the equation to be solve and the method to solve it. 
	\begin{equation}
		\px{u}{t} + \px{f}{x} = 0
	\end{equation}
	If we denote a variable in the reference sub-domain as $\hu_i$, then $\hu_i = J_iu_i$, and for each sub-domain a discontinuous polynomial may be defined as:
	\begin{equation}
		\hud_i = \sum^p_{j=0} \hud_i(\xi_j)l_j(\xi)
	\end{equation}
	where $p$ is the order of the polynomial fit and $l_j(\xi)$ is the $j^{\mathrm{th}}$ Lagrange polynomial.
	\begin{equation}
		l_j(\xi) = \prod_{k=0}^p\bigg(\frac{\xi-\xi_k}{\xi_j-\xi_k}\bigg)(1-\delta_{jk})
	\end{equation}
	The flux polynomial may be similarly defined and subsequently interpolated to the boundary of the element. Having placed boundary nodes on adjacent elements such that they are collocated, a common interface flux may be calculated. The aim of this is to correct the discontinuous polynomial such that a $C^0$ continuous solution is formed. The common interface flux can be calculated by a number of methods, however when solving for hyperbolic equations it is important that the upwind direction is considered in the calculation, and a Riemann problem should be approximately solved. More details are available from Toro~\cite{Toro2009}. With the common interface flux at the left and right interfaces defined($\hfI_l$ and $\hfI_r$ respectively) the flux correction can be applied. This is performed using a correction flux defined such that:
	\begin{equation}\label{eq:flux_cor}
		\hfC(\xi) = (\hfI_l-\hfd_l)h_l(\xi) + (\hfI_r-\hfd_r)h_r(\xi)
	\end{equation} 
	where $h_l$ and $h_r$ are the left and right correction functions. It is obvious that they have boundary conditions: $h_l(-1)=h_r(1) = 1$ and $h_l(1) = h_r(1) = 0$. Further information on the definition of the correction function may be found in \cite{Vincent2010,Vincent2015,Trojak2018}. 
	Upon differentiating Eq.(\ref{eq:flux_cor}) the transformed conservation equation may be written as:
	\begin{align}
		\px{\hu_j}{t} &= -\px{\hf_j}{\xi} \\
					&= -\sum^p_{j=0}\hfd_{i,j}\dx{l_j(\xi)}{\xi} - (\hfI_{i,l}-\hfd_{i,l})\dx{h_l(\xi)}{\xi} - (\hfI_{i,r}-\hfd_{i,r})\dx{h_r(\xi)}{\xi}
	\end{align}
	With the equation semi-discretised, the solution may now be advanced in time via a suitably selected temporal integration method.	
	
\section{Analytical Methods}
	We begin the investigation into the temporal stabilisation of FR via a theoretical numerical analysis. Two analytical techniques will be presented so that the character of FR under the effect of filtering may be assessed. These are von Neumann analysis and an accompanying error analysis. We start by introducing the von Neumann analysis which will be performed in a manner similar to that introduced by Lele~\cite{Lele1992}, Huynh~\cite{Huynh2007}, and Vincent~\etal~\cite{Vincent2011}. The error study aims to extend the work of Hesthaven~\etal~\cite{Hesthaven2007} and Asthana~\etal~\cite{Asthana2017}, with an extension made to show the fully-discretised temporal behaviour of harmonic solutions.
	 
\subsection{Von Neumann Analysis}
	  The von Neumann analysis present will be focused on the linear advection equation.
	\begin{equation}\label{eq:1dadvection}
		\frac{\partial u}{\partial t} + a\frac{\partial u}{\partial x} = 0
	\end{equation}
	A matrix form of this equation consistent with FR operators was presented by Huynh~\cite{Huynh2007} and Vincent~\etal~\cite{Vincent2011} and we will follow this, but with the added generalisations for arbitrary meshes as in Trojak~\etal~\cite{Trojak2017}, which leads to the following form:
	\begin{align}\label{eq:DiscreteFR2}
		\frac{\partial \bar{\mathbf{u}}_j}{\partial t} &= -J_{j+1}^{-1}\mathbf{C}_{+1}\hat{\mathbf{u}}_{j+1} -J_j^{-1} \mathbf{C}_0 \hat{\mathbf{u}}_j - J_{j-1}^{-1} \mathbf{C}_{-1}\hat{\mathbf{u}}_{j-1} \\		
		\mathbf{C}_{+1} &= (1-\alpha)\mathbf{h_r}^T\mathbf{l_l} \\
		\mathbf{C}_0 &= \mathbf{D} - \alpha\mathbf{h_l}^T\mathbf{l_l} - (1-\alpha)\mathbf{h_r}^T\mathbf{l_r}\\		
		\mathbf{C}_{-1} &= \alpha\mathbf{h_l}^T\mathbf{l_r}
	\end{align}
	where $\mathbf{D}$ is the differentiation matrix at the solution points, $\mathbf{h_l}$ and $\mathbf{h_r}$ are the gradients of the left and right correction functions evaluated at the solution points, and $\mathbf{l_l}$ and $\mathbf{l_r}$ are interpolation vectors to the left and right interfaces. The spectral response of the system is then obtained by applying a trial Bloch wave solution and vectorising the result to give: 
	\begin{align}
		u(x,t) &= v e^{i(kx - \omega t)} \label{eq:bloch1d2} \\
		\mathbf{u}_j &= \mathbf{v}_j e^{ik(\frac{\xi+1}{2}\delta_{j+1}+x_j-c(k)t)}\label{eq:bloch1dd2}
	\end{align}
	where $\delta_j=x_j-x_{j-1}$ and for wavenumber $k$. Substituting this into Eq.(\ref{eq:DiscreteFR2}), FR may be cast into the semi-discrete form: 
	\begin{equation}\label{eq:semi_discrete}
		\frac{\partial \bar{\mathbf{u}}_j}{\partial t} = \mathbf{Q} \bar{\mathbf{u}}_j
	\end{equation}
	For the case of generalised interfaces it follows that $\mathbf{Q} = -J_{j+1}^{-1}\mathbf{C}_{+1}e^{ik\delta_{j+1}} - J_j^{-1}\mathbf{C}_0 - J_{j-1}^{-1}\mathbf{C}_{-1}e^{-ik\delta_{j}}$. This may then be transformed into the fully discrete form by application of a suitable temporal integration method, thereby giving what is often called the update equation:
	\begin{align}
			\bar{\mathbf{u}}_j(t + \tau) &= \mathbf{R}(\mathbf{Q})\bar{\mathbf{u}}_j(t)\label{eq:RK_gen}\\
			\bar{\mathbf{u}}_j(n\tau) &= \mathbf{R}^n(\mathbf{Q})\bar{\mathbf{u}}_j(0)\label{eq:RK_gen_n} \\ 
			\mathbf{R}_{33} &= \mathbf{I} + \frac{(\tau\mathbf{Q})^1}{1!} + \frac{(\tau\mathbf{Q})^2}{2!} + \frac{(\tau\mathbf{Q})^3}{3!} \label{eq:RK_33}
	\end{align}
	where $\mathbf{R}$ is the update matrix and is a function of $\mathbf{Q}$ that is dependent on the temporal integration used. In a method similar to that used on the semi-discrete form, the Bloch wave solution can be applied and discretised, leading to the result:
	\begin{align}
		\underbrace{e^{-ik(c-1)\tau}}_{\lambda}\mathbf{v} &= e^{ik\tau}\mathbf{Rv}\label{eq:full_discrete} \\
		e^{-ik(c-1)n\tau}\mathbf{v} &= e^{ikn\tau}\mathbf{R}^n\mathbf{v}\label{eq:full_discrete_n}
	\end{align}
	The result of Eq.(\ref{eq:full_discrete}) was first reported by Asthana and Jameson~\cite{Asthana2014}, but the accompanying analysis was limited, and was followed up by Vermeire and Vincent~\cite{Vermeire2017a} where the application of FR to implicit LES was considered. Equation~(\ref{eq:full_discrete}) can be seen to be an eigenvalue problem of the update matrix, with $\mathbf{v}$ being the eigenmodes of the fully discretised scheme. The eigenvalues of the system are then temporally shifted and scaled solutions, with the fully discretised wave speed recovered as:
	\begin{align}
		c(k;\tau) &= \frac{i\log{(\lambda)}}{k\tau} + 1 \label{eq:temporal_c}\\
		\omega^{\prime} &= c(k;\tau)k \label{eq:temporal_omega}
	\end{align}
	This can provide useful insight into FR, as previously the semi-discreted form led to a disconnect between the numerically realised and analytical dispersion/disspation in all but highly temporally resolved waves. In contrast, Eq.~(\ref{eq:temporal_c}) gives a more complete picture. A further result is that scheme stability can then be understood to be when $\Im{(c(k;\tau))} < 0 ~ \forall ~ k~\in~\mathbb{R}$, which is equivalent to $|\lambda(k;\tau)| \leqslant 1~\forall ~ k ~ \in ~ \mathbb{R}$. This will give the same results as von Neumann theorem for stability~\cite{Isaacson1994}, but will give added information about the mechanism driving the scheme unstable.
	
\subsection{Error Analysis}\label{sec:conv_method}
	Von Neumann analysis can be insightful when trying to understand how a method will react when applied to problem. However, this only gives a snapshot of the behaviour, and it can also be of use to understand how a solution will develop as the dispersion/dissipation transfer function is repeatedly applied and integration error has iterations in which to grow. We also consider, therfore, how the error developes with time, which was initially investigated for Nodal DG via FR by Asthana~\etal~\cite{Asthana2017}. In the semi-discrete sense this may be can be done by utilising exponential integration of Eq.(\ref{eq:semi_discrete}), as $\mathbf{Q}$ forms a linear spatial differentiation operator, see Pope~\cite{Pope1963}. The first step, following the method of~\cite{Asthana2017}, is to diagonalise the operator matrix: 
	\begin{equation}\label{eq:diagonalisation}
		\mathbf{Q} = \mathbf{W}_{Q}\mathbf{\Gamma}_{Q}\mathbf{W}^{-1}_{Q} = ik\mathbf{W}_{Q}\mathbf{\Lambda}_{Q}\mathbf{W}^{-1}_{Q}
	\end{equation}  
	where $\mathbf{W}$ is a matrix of eigenvectors, $\mathbf{\Gamma}$ is a diagonal matrix of eigenvalues and $\mathbf{\Lambda}$ is a diagonal matrix of complex convective velocities. From here the temporal integration can be performed:
	\begin{equation}\label{eq:t_int_1}
		\mathbf{u}_j^{\delta}(t) = \exp{(ct\mathbf{Q})}\mathbf{u}(0) = \mathbf{W}_{Q}\exp{(ikct\mathbf{\Lambda}_Q)}\mathbf{W}^{-1}_{Q}\mathbf{u}_j(0)	
	\end{equation}
	where the definition of the exponential has been used to simplify the form by considering $\mathbf{W}^{-1}\mathbf{W}=\mathbf{I}$.
	By applying the initial condition as a sum of the eigenvectors, again following Asthana~\cite{Asthana2017}:
	\begin{equation}\label{eq:initial_condition}
		\mathbf{u}_j(0) = \exp{(ikx_j)}\mathbf{W}_Q\pmb{\beta} = \exp{(ikx_j)}\sum_{n=0}^{p} \mathbf{w}_{Q,n}\beta_n = \exp{\big(ikJ_j(\pmb{\xi}+1)\big)}\exp{(ikx_j)}
	\end{equation}
	where $\pmb{\beta}$ is an array of weights associated with the reconstruction of the initial solution from the eigenvectors $\mathbf{w}_Q$. Substituting this into Eq.(\ref{eq:t_int_1}) leads to:
	\begin{equation}
		\mathbf{u}_j^{\delta}(t) = \mathbf{W}_{Q}\exp{(ikct\mathbf{\Lambda}_Q)}\exp{(ikx_j)}\pmb{\beta} = \exp{(ikx_j)}\sum_{n=0}^{p}\exp{(ikct\lambda_n)}\beta_n\mathbf{w}_{Q,n}
	\end{equation}
	Lastly, to understand the rate convergence, the error is calculated via analytically evolving the initial condition such that:
	\begin{equation}
		\mathbf{u}_j(t) = \exp{(ik(x_j - ct))}\sum^p_{n=0}\beta_n\mathbf{w}_{Q,n}
	\end{equation}
	and hence the semi-discretised error is:
	\begin{equation}
		\mathbf{e}_j(t) = \mathbf{u}_j^{\delta}(t) - \mathbf{u}_j(t) = \exp{(ik(x_j - ct))}\sum^p_{n=0}\Big(\exp{\big(ikct(\lambda_n+1)\big)}-1\Big)\beta_n\mathbf{w}_{Q,n}
	\end{equation}	
	From here the it is necessary to calculate the eigenvectors and values and then decompose the solution to calculate $\beta_i$, from Eq.(\ref{eq:initial_condition}).
			
	If instead the diagonalised form of the discrete FR operator, Eq.(\ref{eq:diagonalisation}), is substituted into the fully discretised version of the update equation, Eq.(\ref{eq:RK_gen_n}), then:
	\begin{align}
		\mathbf{u}^{\delta}_j(n\tau) &= \Bigg(\sum_{m=0}^{r}{\frac{(ik\tau)^m}{m!}\mathbf{W}\mathbf{\Lambda}^m \mathbf{W}^{-1} } \Bigg)^n\mathbf{u}_j(0) \\
		\mathbf{u}^{\delta}_j(n\tau) &= \Bigg(\sum_{m=0}^{r}{\frac{(ik\tau)^m}{m!}\mathbf{W}\mathbf{\Lambda}^m \mathbf{W}^{-1} } \Bigg)^n\mathbf{W}\pmb{\beta}\exp{(ikx_j)}
	\end{align}	
	And therefore, the fully discrete error can be derived as:
	\begin{equation}\label{eq:FR_error_ts}
		\mathbf{e}_{j,r}(n\tau) = \mathbf{u}^{\delta}_j(n\tau) - \mathbf{u}_j(n\tau) = \Bigg(\Bigg(\sum_{m=0}^{r}{\frac{(ik\tau)^m}{m!}\mathbf{W}\mathbf{\Lambda}^m \mathbf{W}^{-1} } \Bigg)^n - \exp{(-ikcn\tau)\mathbf{I}} \Bigg) \mathbf{W}\pmb{\beta}\exp{(ikx_j)}
	\end{equation}
	where $r$ is the number of RK temporal integration sub-steps and $\mathbf{I}$ is the identity matrix. What this shows is how the truncation of the temporal integration introduces error, and, interestingly that $\lim{(\mathbf{e}_{j,r})}_{r\rightarrow\infty} = \mathbf{e}_j$, by remembering the recursive definition of $\exp{(z)}$. An attempt can be made to form a temporal integration method that reduces the error via the introduction of factors at each RK step, however it was found that this method reduces the CFL limit by an order of magnitude due to the weights varying for each wavenumber, and does not seem to be feasible in the general case.

\section{Filtering}
	Filtering can be used in LES to suppress numerical errors by excluding erroneous high wavenumbers, such that sub-grid scale convergence is separated from numerical convergence, see Ghousal~\cite{Ghosal1996}. This allows explicit LES to converge upon the spatially filtered NSE. This does, however, require the domain-level setting of filter width. This is not directly congruous with FR, where operations are defined by matrices that are spatially invariant due to the use of a reference domain. Therefore, it is proposed that reference domain filtering is used such that the computational benefit of spatially invariant operators may still be realised, whilst still using filtering to modify the underlying scheme's behaviour. 
	
	Filtering was recently investigated by Asthana~\etal~\cite{Asthana2015a}, where it was used with a globally defined filter width to suppress aliasing errors when non-linear equations are considered. It was proven that filtering of the solution in FR is equivalent to added dissipation and super-dissipation. It is therefore thought that local reference domain filtering may allow control over the linear advection dissipation. 
	
	We begin by defining the filter to be studied and the possible means of application. There are a plethora of possible filters, with a list of important ones given by Pope~\cite{Pope2010}, for which the inverse Fourier transform can be used to calculate the spatial filter kernel, $s(r)=\mathcal{F}^{-1}\{S(k)\}$. We will primarily focus on the Gaussian filter kernel, as this is a canonical example which can also be presented in a $C^p$ continuous manner which is believed may make the result more stable. The filter matrix can be defined as:
	\begin{equation}\label{eq:gaussian_filter}
		\mathbf{S}_{i,j} = \sqrt{\frac{6}{\pi\sigma^2}}\exp{\bigg(-\frac{6|\pmb{\xi}_i-\pmb{\xi}_j|^2}{\sigma^2}\bigg)}
	\end{equation}	 
	where the filter width is $\sigma$ and the transformed solution point coordinates are $\pmb{\xi}$. To ensure that the filter does not act as a net source or sink the integrated effect of each solution point is used to renormalise the effect. The filter can be applied in three ways, with varying affect. The first approach is to apply the filter to the entire spatial scheme, as in:
	\begin{equation}\label{eq:Q_filter}
		\frac{\partial \bar{\mathbf{u}}_j}{\partial t} = \mathbf{SQ} \bar{\mathbf{u}}_j
	\end{equation}
	The filter may also be just applied to the differentiation operator, defining $\mathbf{C}_0$ as:
	\begin{equation}\label{eq:D_filter}
		\mathbf{C}_0 = \mathbf{SD} - \alpha\mathbf{h_l}^T\mathbf{l_l} - (1-\alpha)\mathbf{h_r}^T\mathbf{l_r}
	\end{equation}
	which would allow the correction function to introduce higher wavenumber content but suppress it from the differentiation operator. Finally, the filter can be exclusively applied to the correction function, defining $\mathbf{C}_0$, $\mathbf{C}_{-1}$ and $\mathbf{C}_{+1}$ as:
	\begin{align}	
		\mathbf{C}_{+1} &= (1-\alpha)\mathbf{Sh_r}^T\mathbf{l_l} \label{eq:C+_filter}\\
		\mathbf{C}_0 &= \mathbf{D} - \alpha\mathbf{Sh_l}^T\mathbf{l_l} - (1-\alpha)\mathbf{Sh_r}^T\mathbf{l_r} \label{eq:C0_filter}\\		
		\mathbf{C}_{-1} &= \alpha\mathbf{Sh_l}^T\mathbf{l_r} \label{eq:C-_filter}
	\end{align}
	which would have the effect of convolving the correction function by the filter therefore preventing the correction function from passing high wavenumber information from adjacent elements. The advantage of the filter approaches proposed in Eq.(\ref{eq:C+_filter}-\ref{eq:C-_filter}) is that the scheme can be potentially filtered in specific areas that such that the ingress of spurious or unwanted modes maybe suppressed.

\section{Results}
	\subsection{Von Neumann Analysis}\label{sec:vn_results}
\begin{figure}
		\centering
		\begin{subfigure}[b]{0.45\linewidth}
			\centering
			\includegraphics[width=\linewidth]{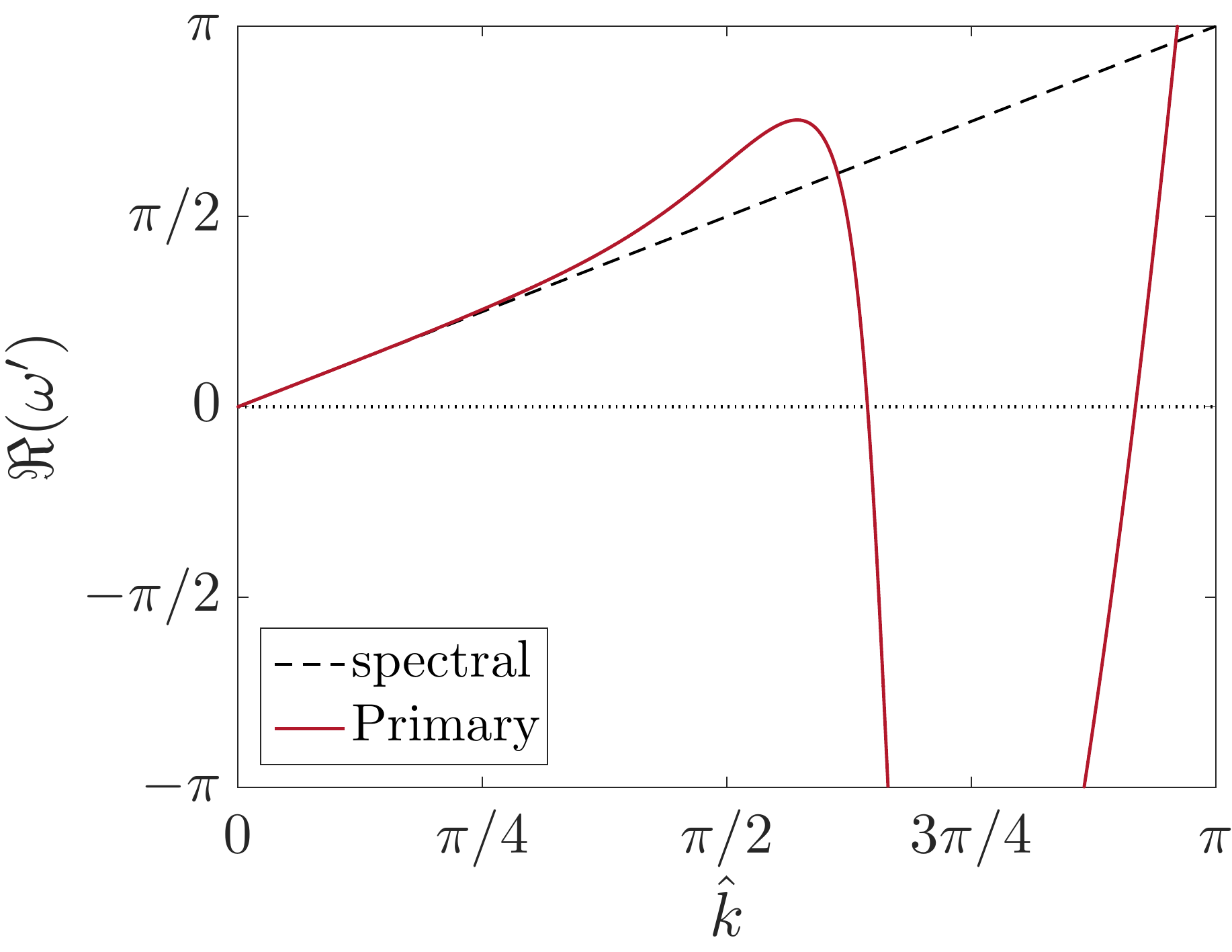}
			\caption{Dispersion\\ \quad \\ \quad }
			\label{fig:FR4_33_real_166}
		\end{subfigure}
		~
		\begin{subfigure}[b]{0.45\linewidth}
			\centering
			\includegraphics[width=\linewidth]{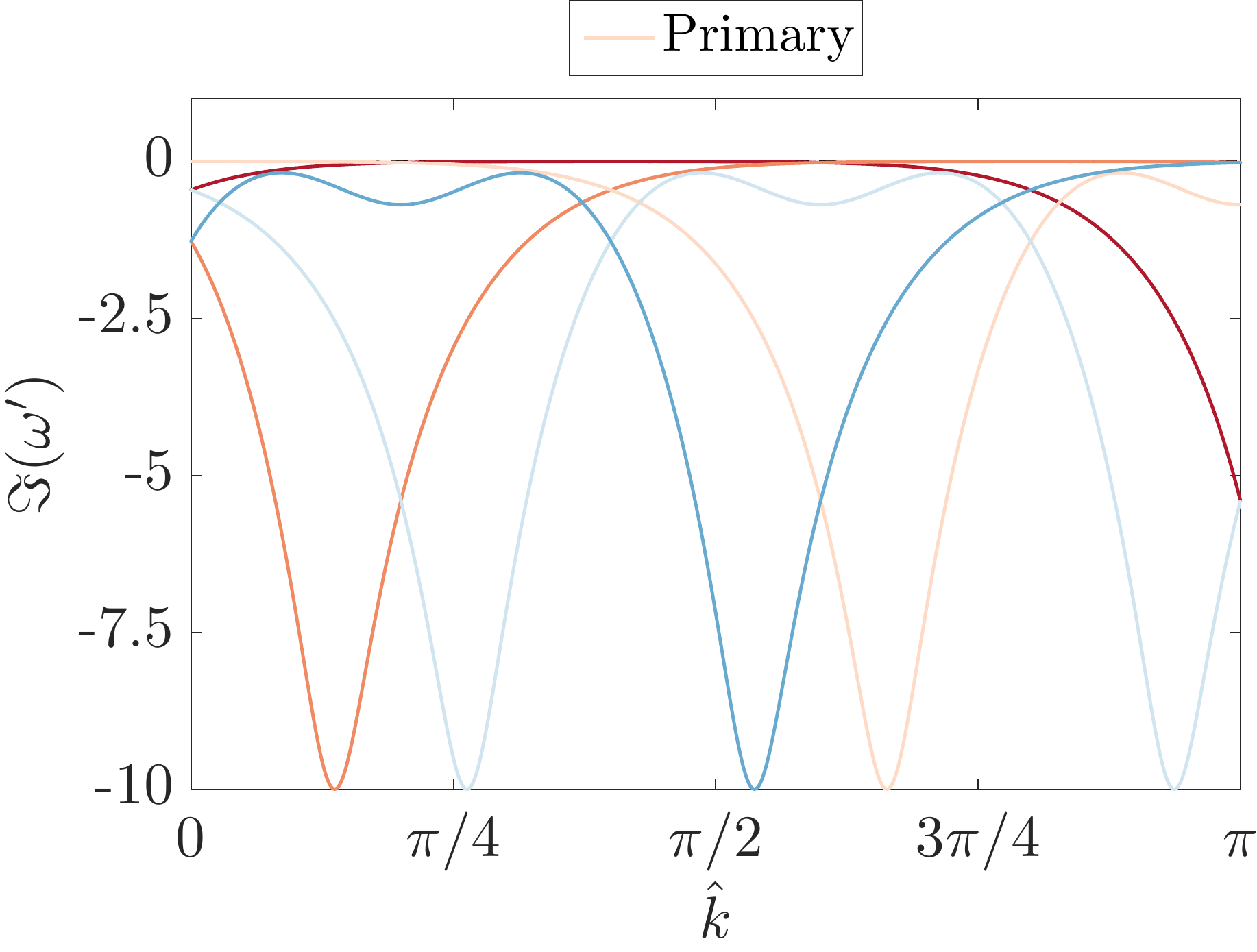}
			\caption{Dissipation\\ \quad \\ \quad }
			\label{fig:FR4_33_imag_166}
		\end{subfigure}
		\caption{Fully discretised analysis of FR using Huynh correction functions, $p=4$ with RK33 time integration, $\tau=0.166$, close to the temporal stability limit.}
		\label{fig:FR4_33_t166}
	\end{figure}

\begin{figure}
		\centering
		\begin{subfigure}[b]{0.45\linewidth}
			\centering
			\includegraphics[width=\linewidth]{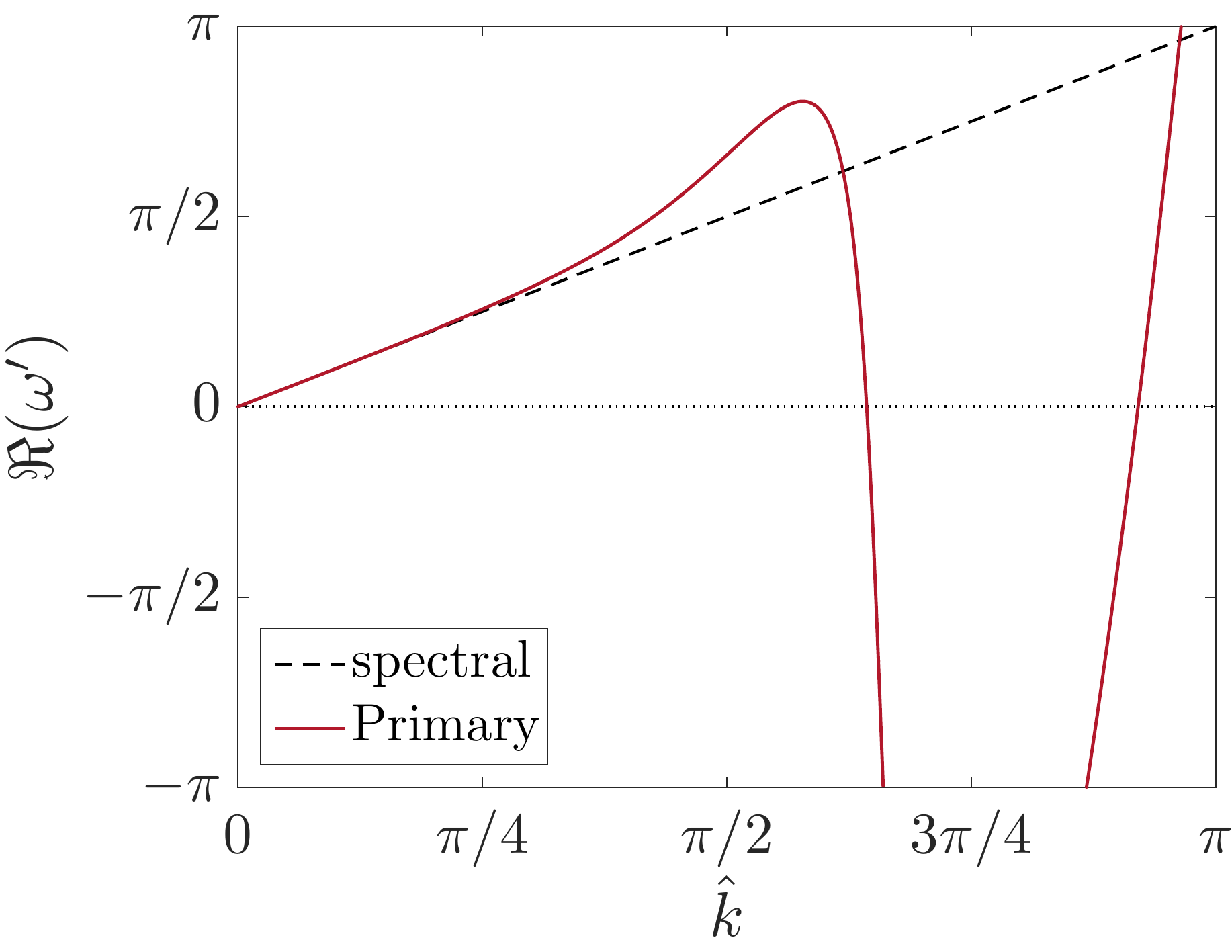}
			\caption{Dispersion\\ \quad \\ \quad }
			\label{fig:FR4_33_real_170}
		\end{subfigure}
		~
		\begin{subfigure}[b]{0.45\linewidth}
			\centering
			\includegraphics[width=\linewidth]{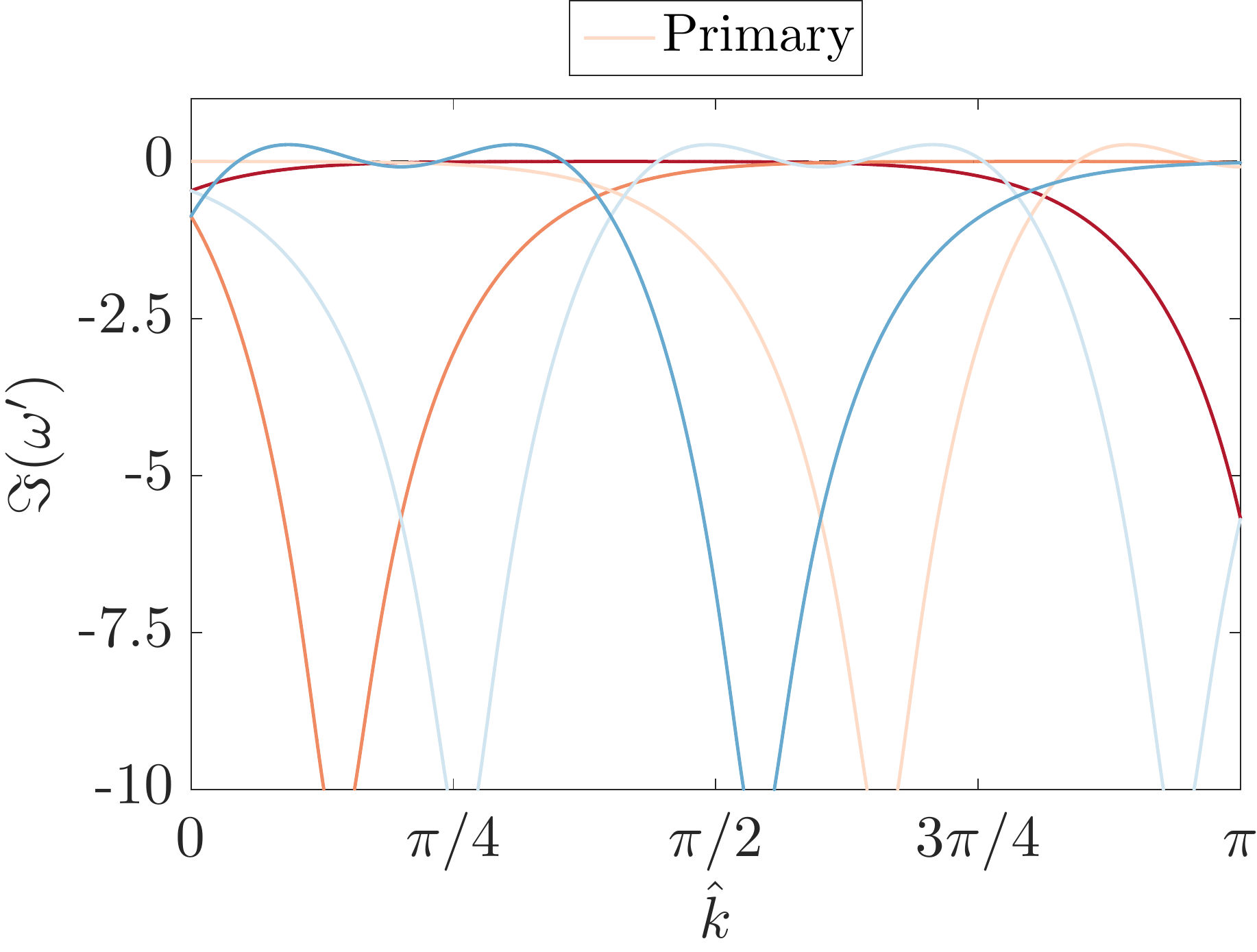}
			\caption{Dissipation\\ \quad \\ \quad }
			\label{fig:FR4_33_imag_170}
		\end{subfigure}
		\caption{Fully discretised analysis of FR using Huynh correction functions, $p=4$ with RK33 time integration, $\tau=0.17$. Showing instability due to exceeding the CFL limit.}
		\label{fig:FR4_33_t170}
	\end{figure}

		The main aim of this paper is to investigate the temporal stability of the scheme and hence the dissipation of the scheme is of primary importance. However, as was observed by Vermiere and Vincent~\cite{Vermeire2017a}, the group velocity of waves near the Nyquist limit becomes very high when near to the CFL limit. Hence, a further result of interest would be the control of this. Initially evaluating the performance of the scheme without filtering, using RK33 temporal integration, it can be seen by comparing Fig.\ref{fig:FR4_33_t166} and Fig.\ref{fig:FR4_33_t170}, with $\tau$ either side of the CFL limit for uniform grids, that, as expected, the temporal instability originates from higher wavenumbers. However, instability is spread to a wider band of wavenumbers via the other eigenmodes. 
				
		However,  for general solutions a multiplicity of modes will form a solution. Therefore, the instability shown Fig.\ref{fig:FR4_33_t170} will likely result in growth of errors and ultimately lead to failure. But, what can be seen is that the demanded harmonic mode only becomes anti-dissipative at high wavenumbers after a region of large dissipation. This makes the suppression of unstable high wavenumber modes via filtering viable, as solution propagation into this high wavenumber region is already bottle-necked by the high dissipation for  $\hat{k} \approx 5\pi/8$. Hence, limiting the effect of filtering on solution accuracy at high wavenumbers.
		
\begin{figure}
		\centering
		\begin{subfigure}[b]{0.45\linewidth}
			\centering
			\includegraphics[width=\linewidth]{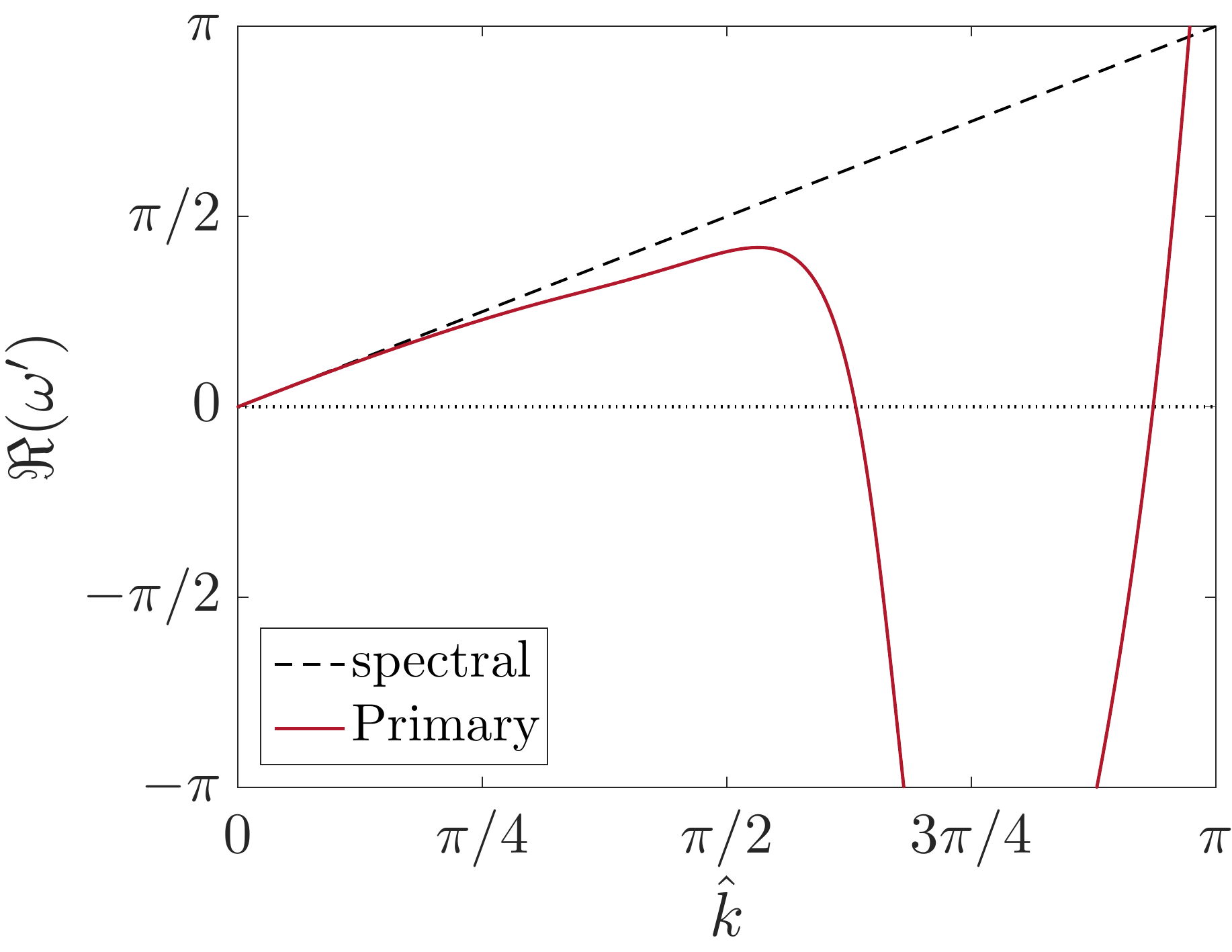}
			\caption{Dispersion\\ \quad \\ \quad }
			\label{fig:FR4_33_real_170_06}
		\end{subfigure}
		~
		\begin{subfigure}[b]{0.45\linewidth}
			\centering
			\includegraphics[width=\linewidth]{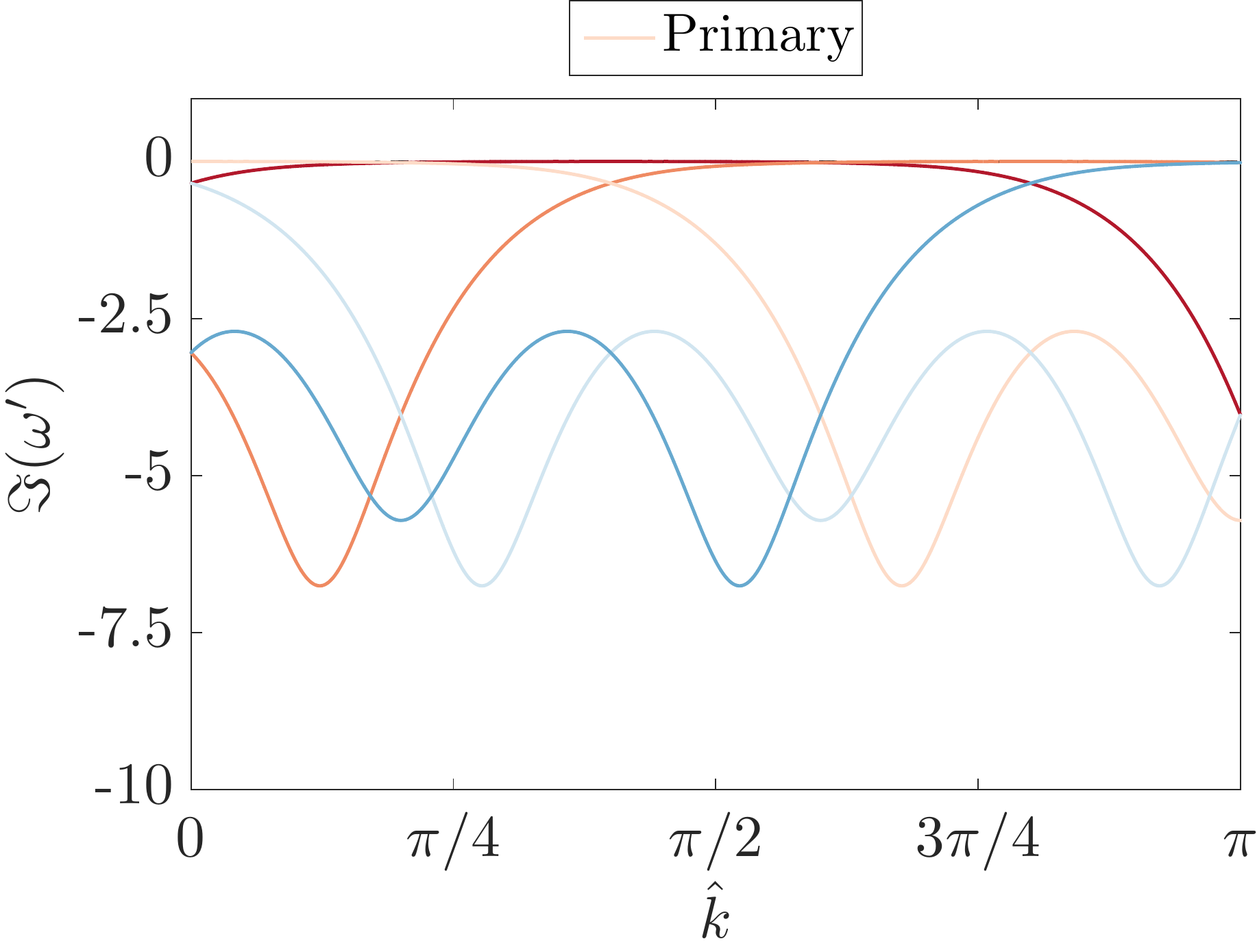}
			\caption{Dissipation\\ \quad \\ \quad }
			\label{fig:FR4_33_imag_170_06}
		\end{subfigure}
		\caption{Fully discretised analysis of FR using Huynh correction functions, $p=4$ with RK33 time integration, $\tau=0.17$. Gaussian spatial filtering is applied as in Eq.(\ref{eq:Q_filter}) with $\sigma=0.6$.}
		\label{fig:FR4_33_t170_06}
	\end{figure}

	Now we apply a Gaussian filter to FR as in Eq.(\ref{eq:gaussian_filter}~\&~\ref{eq:Q_filter}).  Analytically, from Fig.~\ref{fig:FR4_33_t170_06}, with $\sigma=0.6$, it can be seen that the scheme has been stabilised for a CFL number that was previously unstable. The impact of the filter is clearly visible, with high wavenumbers in primary harmonic showing added dissipation. Yet, the region of wavenumbers that experience high dissipation has now had that comparatively reduced. The Nyquist diagrams in Figs.\ref{fig:FR4_33_t170}~\&~\ref{fig:FR4_33_t170_06} are useful in visualising what has happened in this case. These show that although peak dissipation has reduced, the overall dissipation has increased, as might be expected.
	
	The trade off that has been made is evident from Fig.~\ref{fig:FR4_33_real_166}~\&~\ref{fig:FR4_33_real_170_06}, where the phase and group velocities for a central band of wavenumbers ($\pi/4<\hat{k}<\pi/2$) has been reduced. This will impact the convective velocity of solution which will become more evident during later error studies. 
	
	\begin{figure}
		\centering
		\begin{subfigure}[b]{0.47\linewidth}
			\centering
			\includegraphics[width=\linewidth]{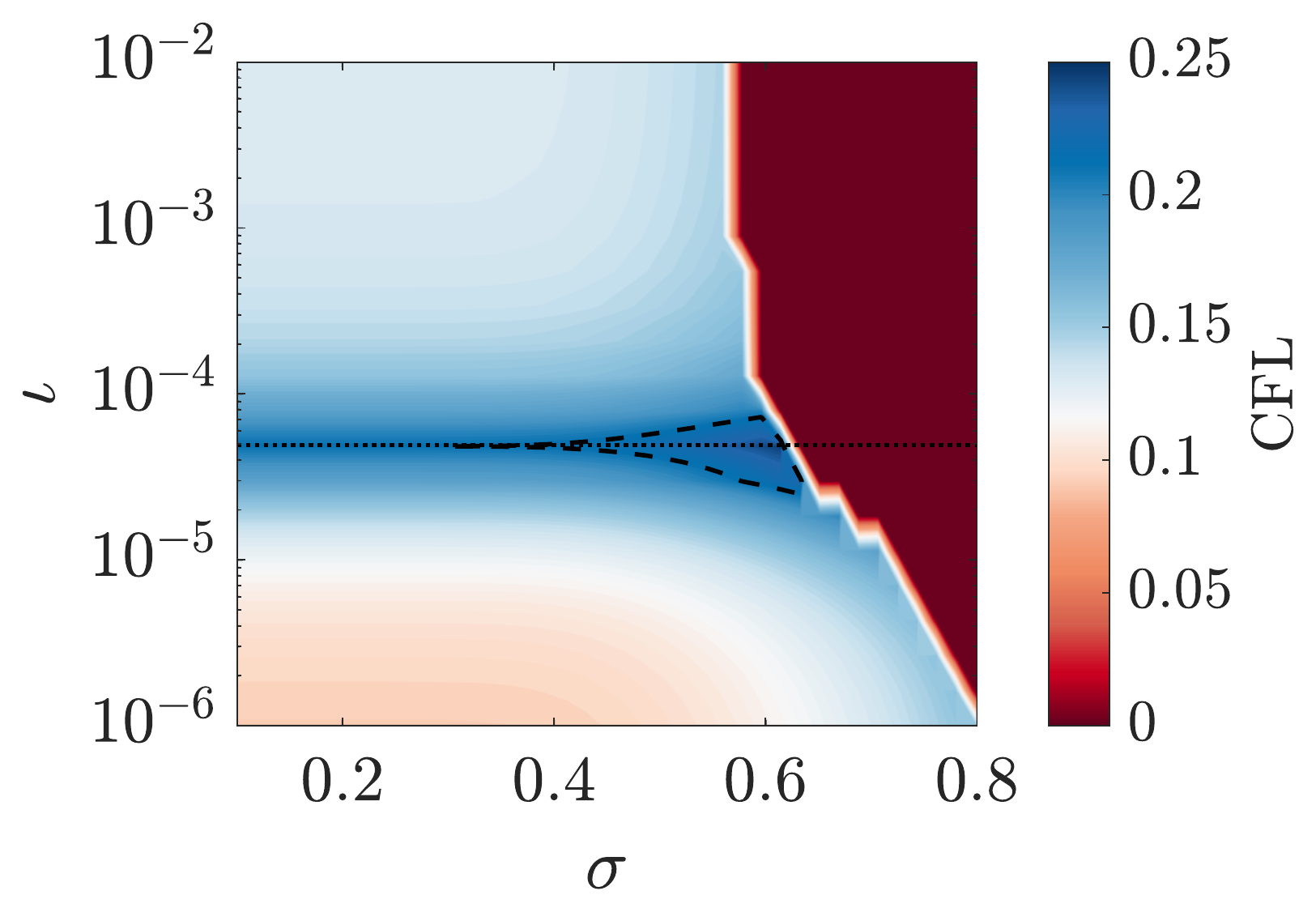}
			\caption{Eq.(\ref{eq:Q_filter})} 
			\label{fig:FR4_RK33_1}
		\end{subfigure}
		~
		\begin{subfigure}[b]{0.47\linewidth}
			\centering
			\includegraphics[width=\linewidth]{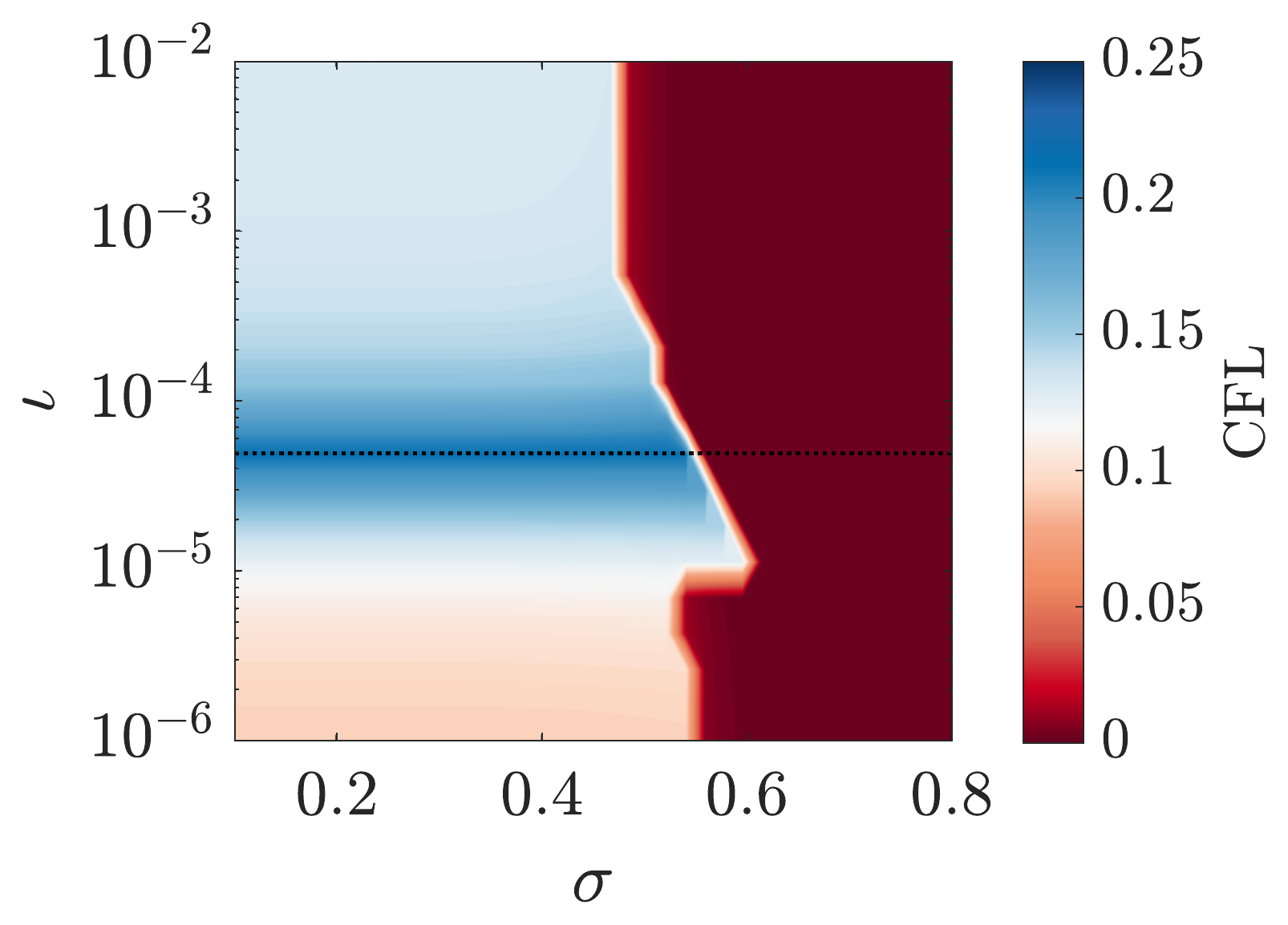}
			\caption{Eq.(\ref{eq:D_filter})} 
			\label{fig:FR4_RK33_2}
		\end{subfigure}
		~
		\begin{subfigure}[b]{0.47\linewidth}
			\centering
			\includegraphics[width=\linewidth]{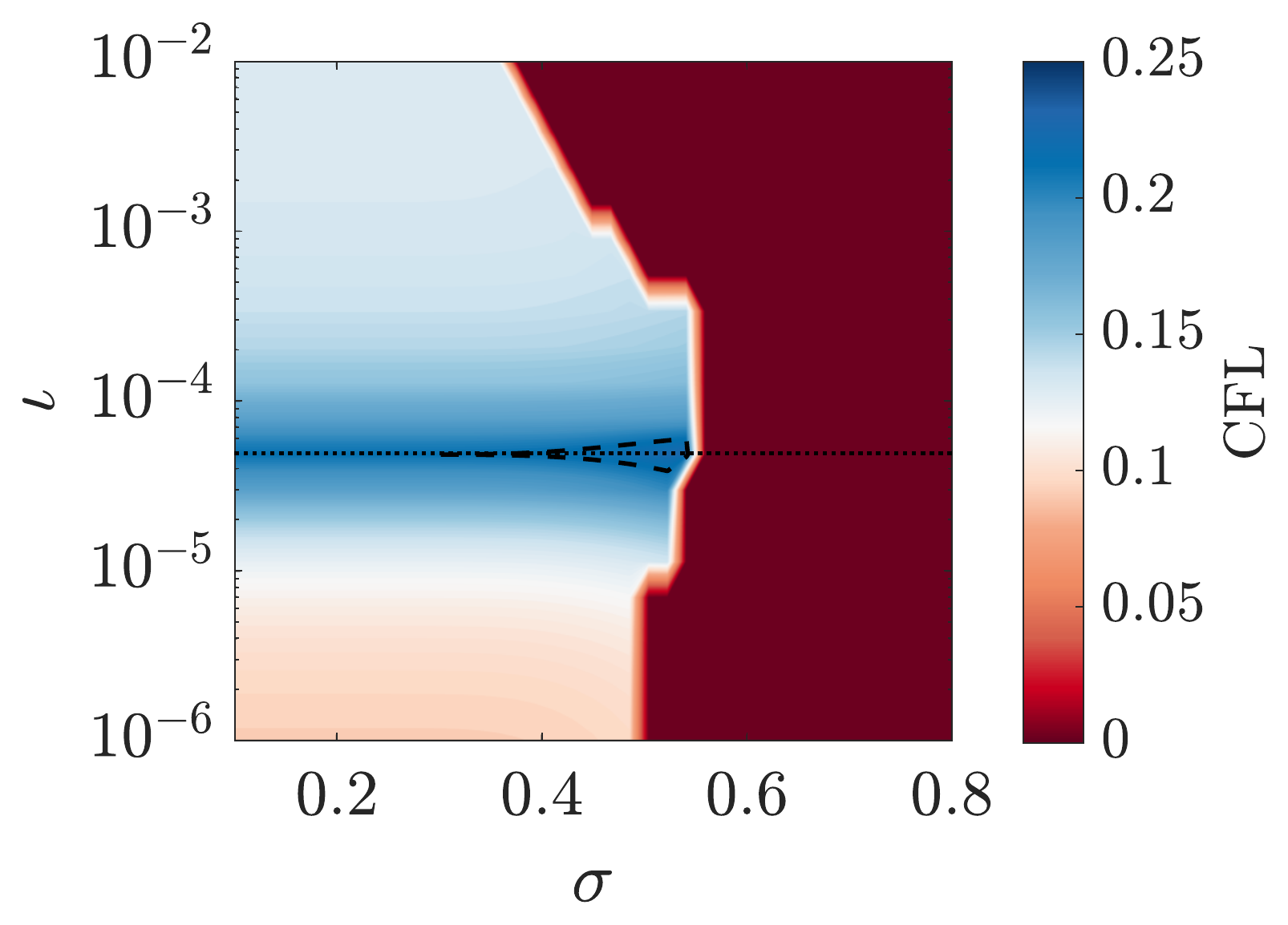}
			\caption{Eq.(\ref{eq:C+_filter}-\ref{eq:C-_filter})} 
			\label{fig:FR4_RK33_3}
		\end{subfigure}		
		\caption{FR, $p=4$, RK33, showing variation of CFL limit with filter width, $\sigma$, and correction function parameter, $\iota$ for various methods of filter application. The dotted line is the $\iota_+$ correction of Vincent~\etal~\cite{Vincent2011} that give peak CFL limit, with the dashed contour at the CFL limit associated with $\iota_+$ ($\mathrm{CFL}=0.2118$).}
		\label{fig:FR4_RK33_cs}
	\end{figure}
	
	This numerical analysis is concluded by considering the von Neumann conditions for ensuring that the operator performs a contraction transformation when temporally discretised, \emph{i.e.} $\rho(\mathbf{R}) \leqslant 1$. The aim of this analysis is to understand the extent to which filtering can improve temporal stability and to understand the maximum stable filter width for the various methods of filter application. For this we will consider only the single parameter energy stable flux reconstruction scheme~\cite{Vincent2010}. It was previously found by Vincent~\etal~\cite{Vincent2011} that there is an optimal correction function that recovers super-convergence, with $\iota=\iota_+$, and gives peak the CFL limit of VJCH correction functions when spatial temporal discretisation is considered.
	
	Through Fig.~\ref{fig:FR4_RK33_cs} the temporal stability of the various methods of filter application are explored. It is clear that filtering the full scheme provides the most temporally stable results, with peak CFL also occurring at $\iota_+$, giving a maximum $25\%$ boost in CFL limit. What is made clear from Fig.~\ref{fig:FR4_RK33_2} is that the application of filtering to solely the differentiation operator has only a limited impact on temporal stability. This is consistent with the order of the differentiation operator being lower than that of the correction and that instability is instigated in the range of wavenumbers that are resolved by the higher order. 
	
	\subsection{Analytic Error Study}\label{sec:filter_conv_results}
\begin{figure}
		\centering
		\begin{subfigure}[b]{0.45\linewidth}
			\centering
			\includegraphics[width=\linewidth]{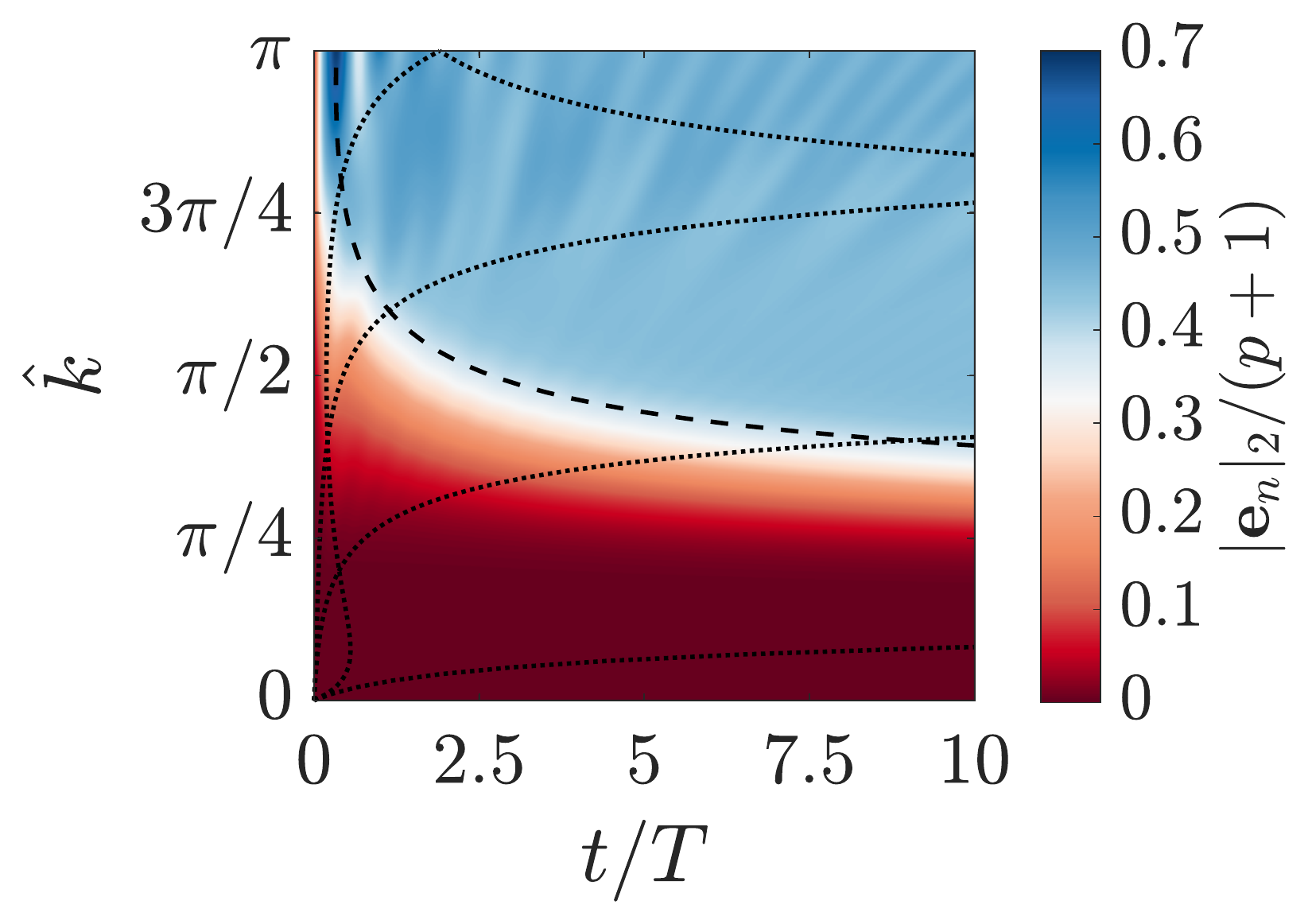}
			\caption{$\sigma = 0$} 
			\label{fig:FR4_conv_00}
		\end{subfigure}
		~
		\begin{subfigure}[b]{0.45\linewidth}
			\centering
			\includegraphics[width=\linewidth]{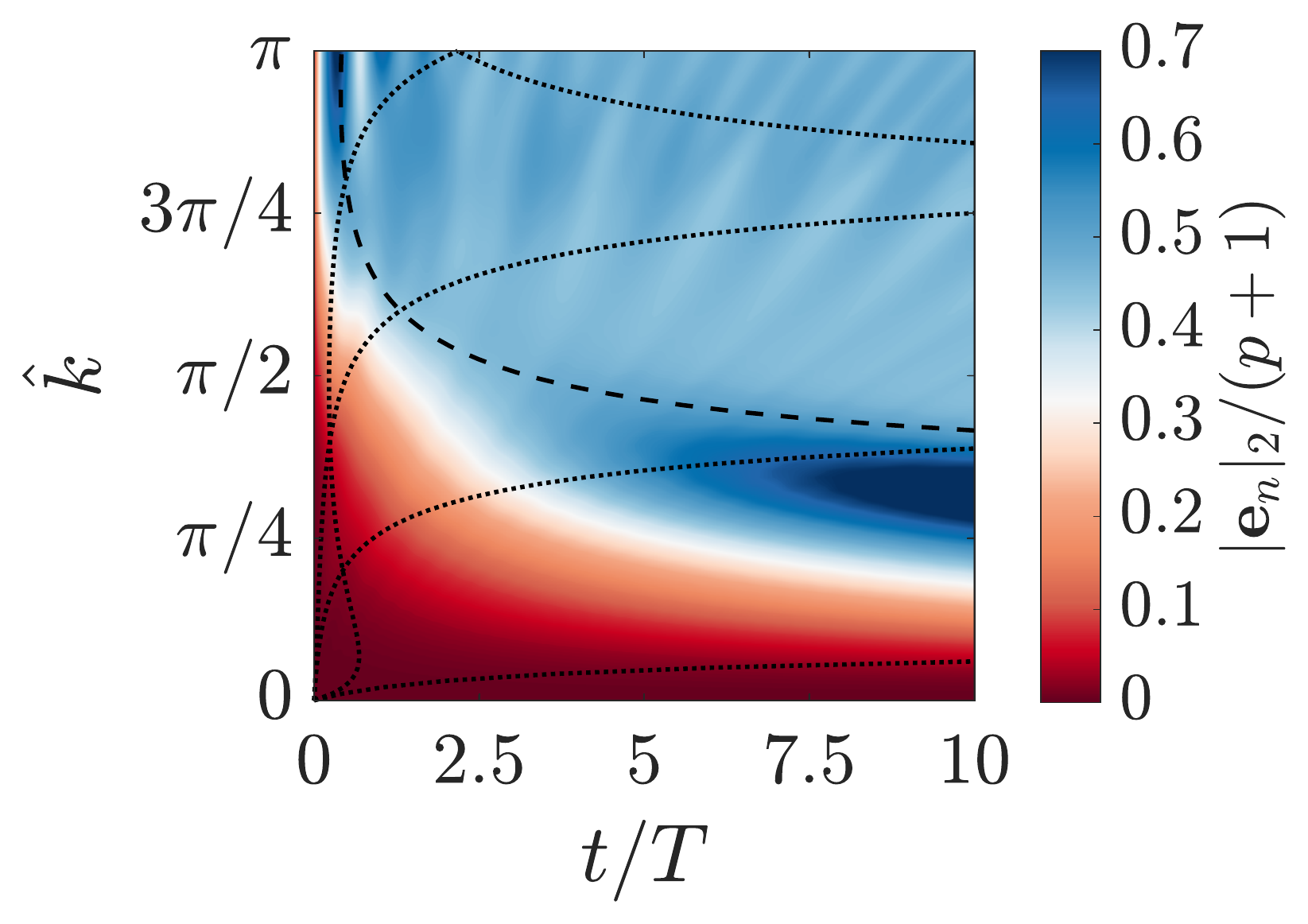}
			\caption{$\sigma= 0.6$} 
			\label{fig:FR4_conv_06}
		\end{subfigure}
		\caption{Analytical error against time and normalised wavenumber. This is for FR with $p=4$ and interface upwinding ($\alpha = 1$) on a uniform mesh. The dashed line is the decay time period for the primary harmonic and the other time periods are shown in dotted lines. }
		\label{fig:FR4_conv}
	\end{figure}

	The effect of a Gaussian filter in the spectral domain is very well understood and it can be seen from Section~\ref{sec:vn_results} that filtering adds dissipation to the high wavenumbers. When considering the temporal stability of the scheme this is very beneficial, however, it is important to consider the effect of filtering on the accuracy of the scheme. Following on from the method introduced in Section~\ref{sec:conv_method}, the semi-discrete error can be calculated, the results of which are shown in Fig.~\ref{fig:FR4_conv}.
	In addition to this analysis, the dissipative half life of the modes is calculated in accordance with Asthana~\etal~\cite{Asthana2017}:
	\begin{equation}
		\tau_{1/2,n} = -\frac{1}{ckJ^{-1}_j\Im{(\lambda_n)}}
	\end{equation}
	where $\tau_{1/2,n}$ is the half life of the $n^{\mathrm{th}}$ mode (the contours of which are shown in Fig.~\ref{fig:FR4_conv}), and $\lambda_n$ is the $n^{\mathrm{th}}$ diagonal entry of $\mathbf{\Lambda}$. Such that $\mathbf{\Lambda} = \bigoplus_{n=0}^p\lambda_n$. From initially studying the unfiltered scheme error, Fig.\ref{fig:FR4_conv_00}, the resolvable wavenuimbers can be broken down into three regions, firstly a region where the half-life of the primary mode is long. In this region waves are well resolved even in the long time integration. Secondly, at slightly higher wavenumbers the primary mode decay time decreases with out a large increase in the half-life of the other modes, in this region the initially well resolved solution will decay uniformly. Lastly at the highest wavenumbers there is a region where the primary mode half-life is low but the half-0lives of the the other modes is large. This leads to the the decay of the primary mode, but slower moving secondary modes do not decay as readily. Hence, stripping error in the error plot. 
	
	Studying the filtered case, Fig.\ref{fig:FR4_conv_06}, much of the same behaviour can be seen, however there is an intermediate band in the range of wavenumbers from $\pi/4<\hat{k}<\pi/2$ where the error seems to become very high. This, by comparison with Fig.\ref{fig:FR4_33_real_170_06}, is the region in which there is dispersion undershoot. Hence in this region the phase and group velocities are lower and the error will increase as the phase difference increases. As the primary mode half-life is broadly unaffected in this region it is expected that at some time later the wave will have locked on again with this analytic solution and the error will decrease. For the semi-discrete case the remaining behaviour at higher wavenumbers looks much the same as the unfiltered case.

	\begin{figure}
		\centering
		\begin{subfigure}[b]{0.45\linewidth}
			\centering
			\includegraphics[width=\linewidth]{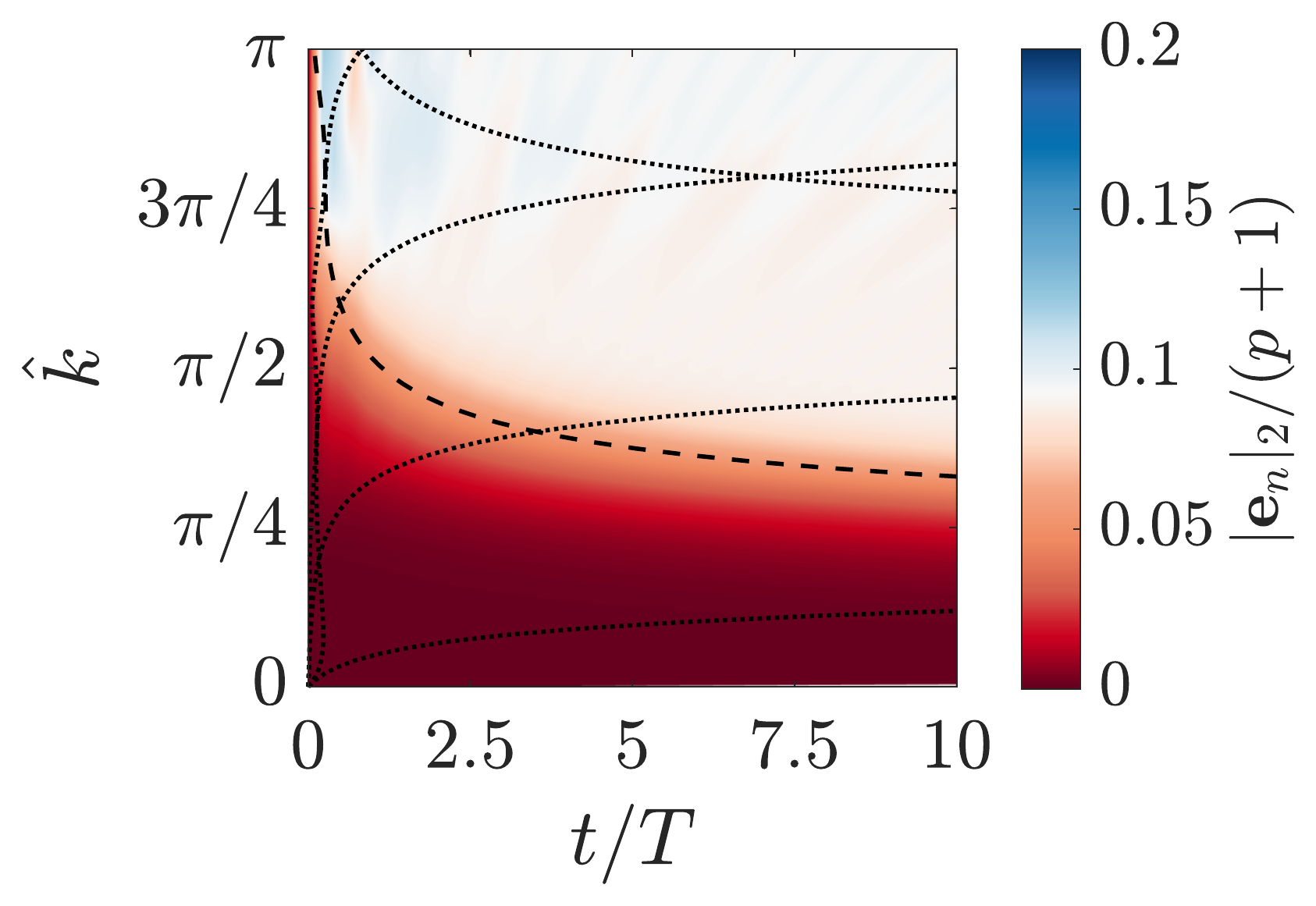}
			\caption{$\sigma = 0,\tau = 0.1$} 
			\label{fig:FR4_conv_ts10_00}
		\end{subfigure}
		~
		\begin{subfigure}[b]{0.45\linewidth}
			\centering
			\includegraphics[width=\linewidth]{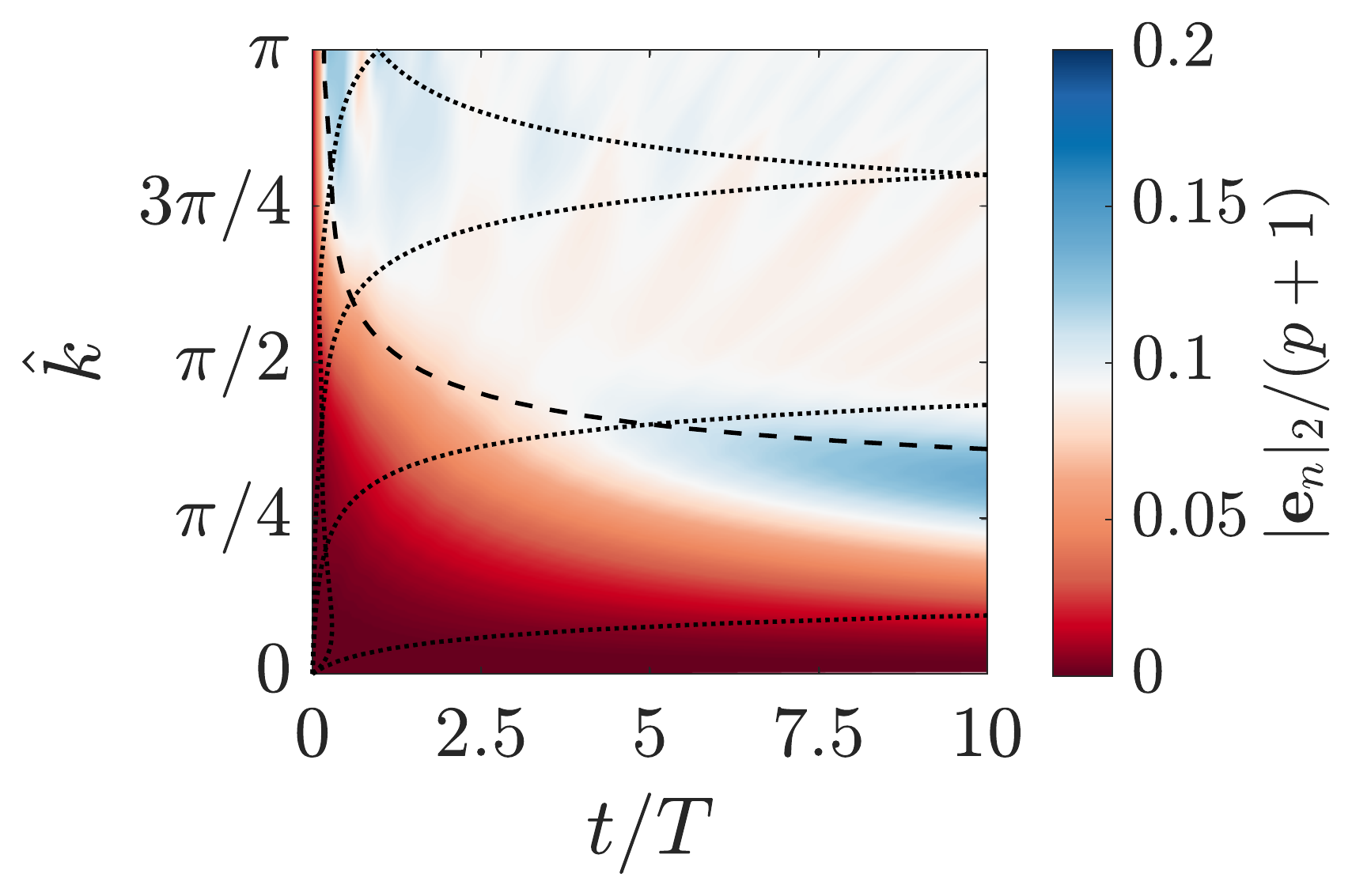}
			\caption{$\sigma= 0.6,\tau = 0.1$} 
			\label{fig:FR4_conv_ts10_06}
		\end{subfigure}
		~
		\begin{subfigure}[b]{0.45\linewidth}
			\centering
			\includegraphics[width=\linewidth]{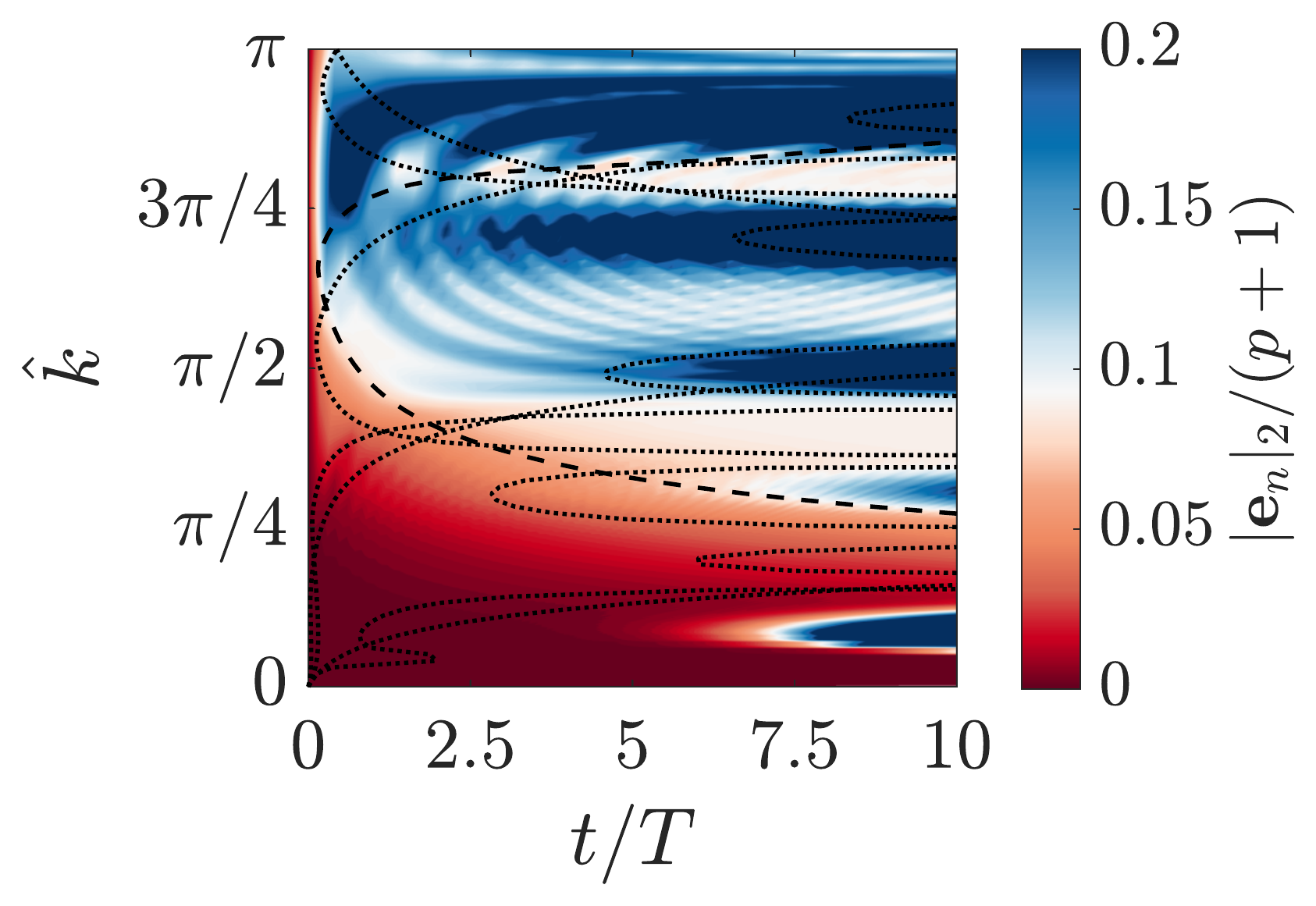}
			\caption{$\sigma = 0,\tau = 0.17$} 
			\label{fig:FR4_conv_ts17_00}
		\end{subfigure}		
		~
		\begin{subfigure}[b]{0.45\linewidth}
			\centering
			\includegraphics[width=\linewidth]{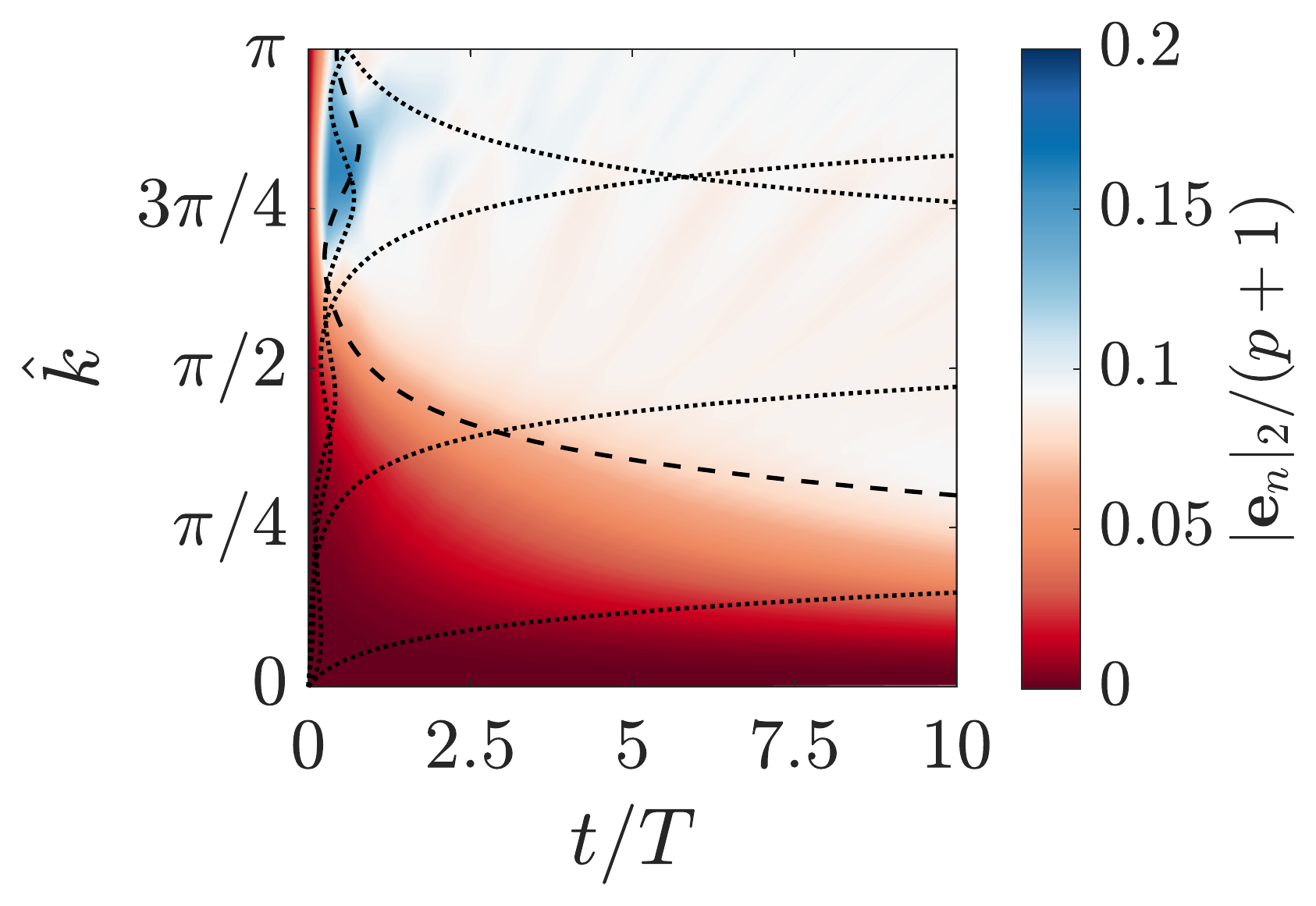}
			\caption{$\sigma = 0.6,\tau = 0.17$} 
			\label{fig:FR4_conv_ts17_06}
		\end{subfigure}	
		\caption{Analytical error against time and normalised wavenumber, with various filters. This is for FR with $p=4$ and interface upwinding ($\alpha = 1$) on a uniform mesh, when temporally and spatially discretised with RK33 temporal integration. Two diffenrent time steps are used, which are stable and unstable for the Huynh correction functions used here. The dashed line is the decay time period for the primary harmonic and the other time periods are shown in dotted lines.}
		\label{fig:FR4_conv_ts}
	\end{figure}
	
	Using Eq.(\ref{eq:FR_error_ts}) the fully discretised error can be explored, along with how discrete temporal integration will effect the error production mechanisms. Comparison of Fig.\ref{fig:FR4_conv_ts10_00} and Fig.\ref{fig:FR4_conv_ts17_00} shows the mechanism of solution divergence when the  CFL limit is exceeded. There are five distinct zones of wavenumbers that become unstable and widen as the time progresses. These regions of instability align with the unstable dissipation that is displayed in Fig.\ref{fig:FR4_33_imag_170}. If the effect of filtering is then considered, as in Fig.\ref{fig:FR4_conv_ts10_06}\&\ref{fig:FR4_conv_ts17_06}, similar behaviour of the error is observed as in the case of exponential filtering. However, at high wavenumbers the region where the dissipation in Fig.\ref{fig:FR4_33_imag_170_06} approaches zero is highlighted as a region where the primary mode half-life is slightly increased. Clearly this is benificial, as it delays the onset of spurious modes. What may also be noted at these higher wavenumbers is that the higher frequency spurious modes are less pronounced when discrete temporal integration is considered, which originates from the discrete temporal integration also behaving some what like a filter. 
		
	\subsection{Numerical Experiments}
\begin{figure}
		\centering
		\begin{subfigure}[b]{0.45\linewidth}
			\centering
			\includegraphics[width=\linewidth]{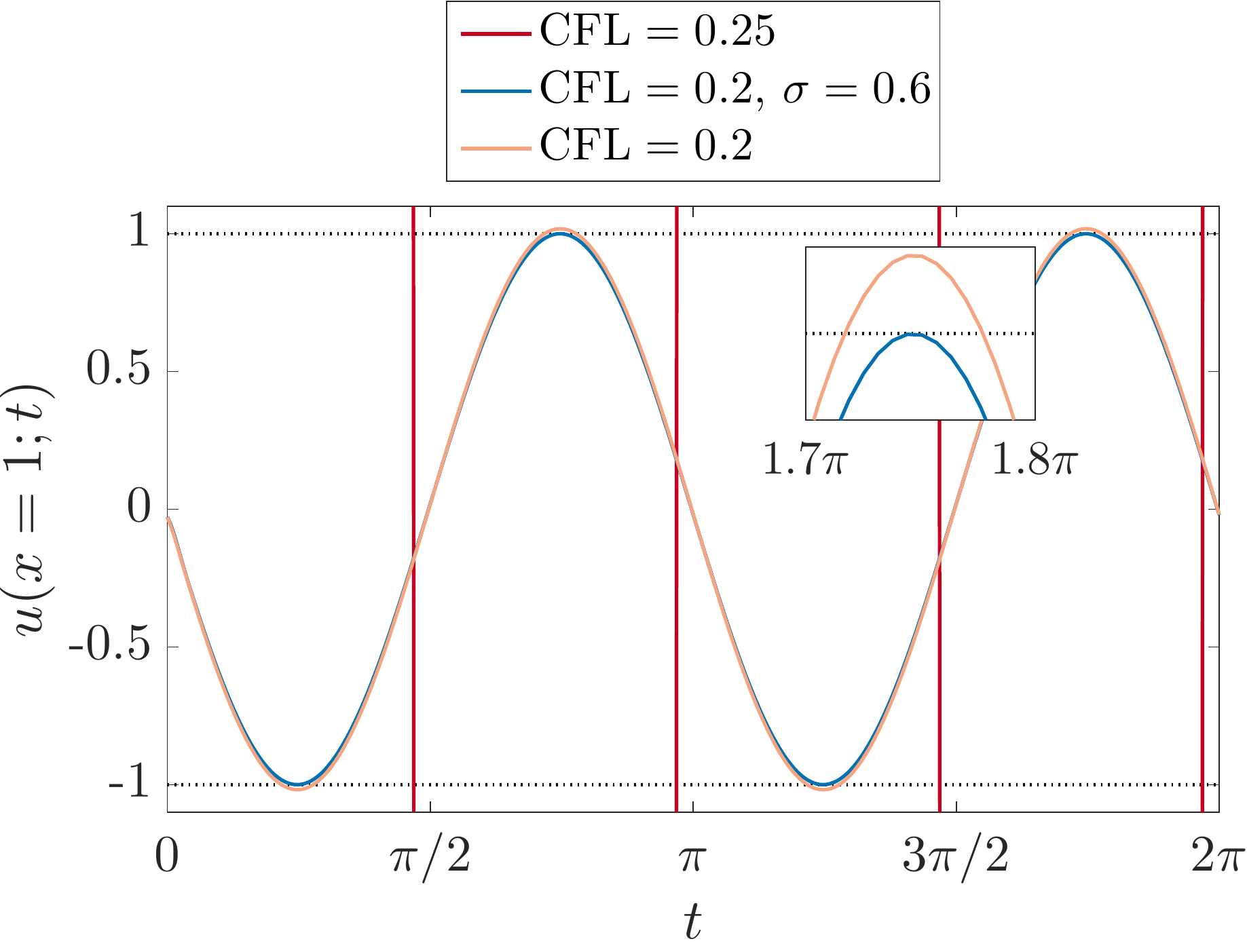}
			\caption{Output in time at $x=1$ with $k=1$.} 
			\label{fig:FR4_33_k2}
		\end{subfigure}
		~
		\begin{subfigure}[b]{0.45\linewidth}
			\centering
			\includegraphics[width=\linewidth]{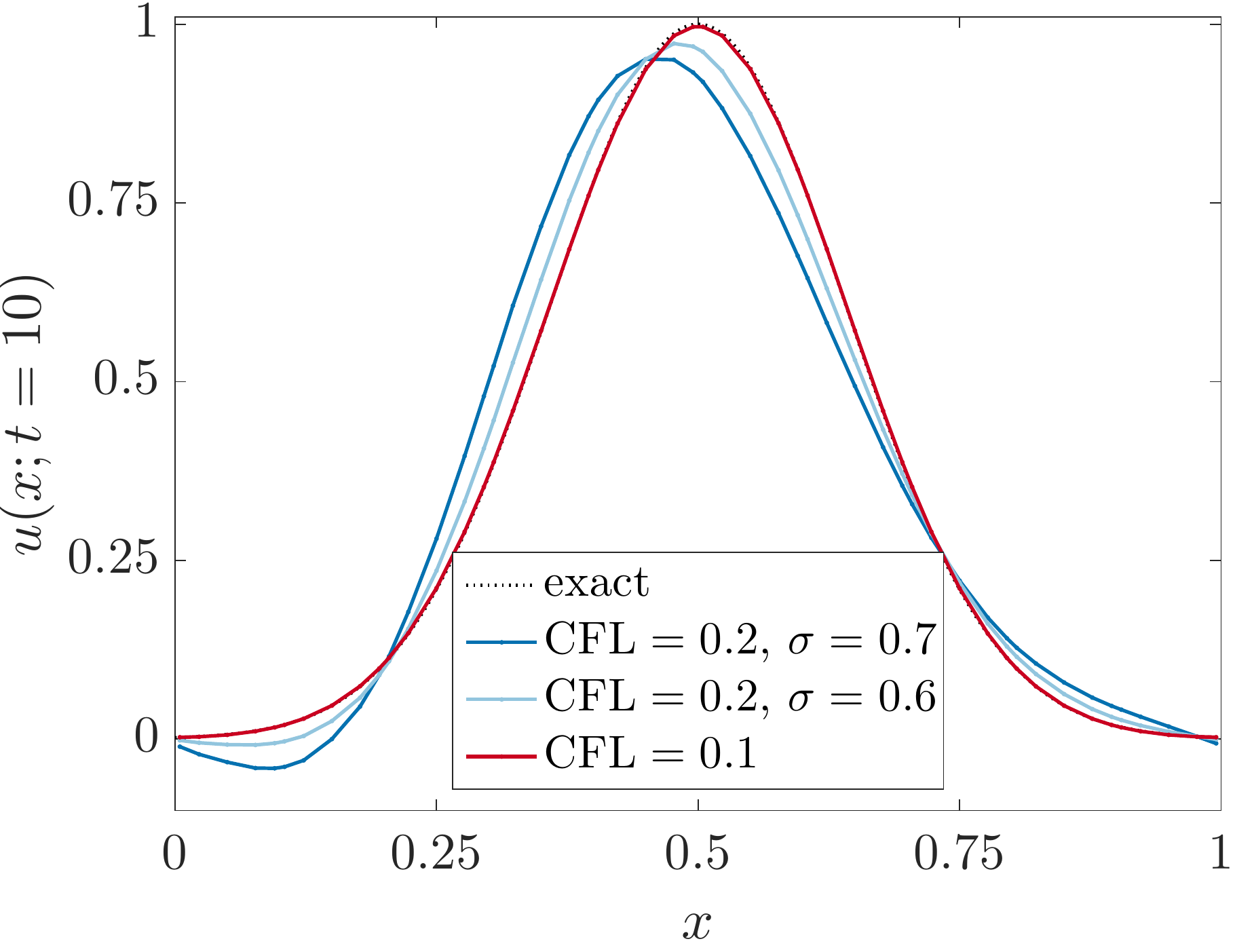}
			\caption{Spatial slice of Gaussian bump, with periodic boundary conditions, after 10 cycles on a unit domain.} 
			\label{fig:FR4_33_bump}
		\end{subfigure}
		~
		\begin{subfigure}[b]{0.45\linewidth}
			\centering
			\includegraphics[width=\linewidth]{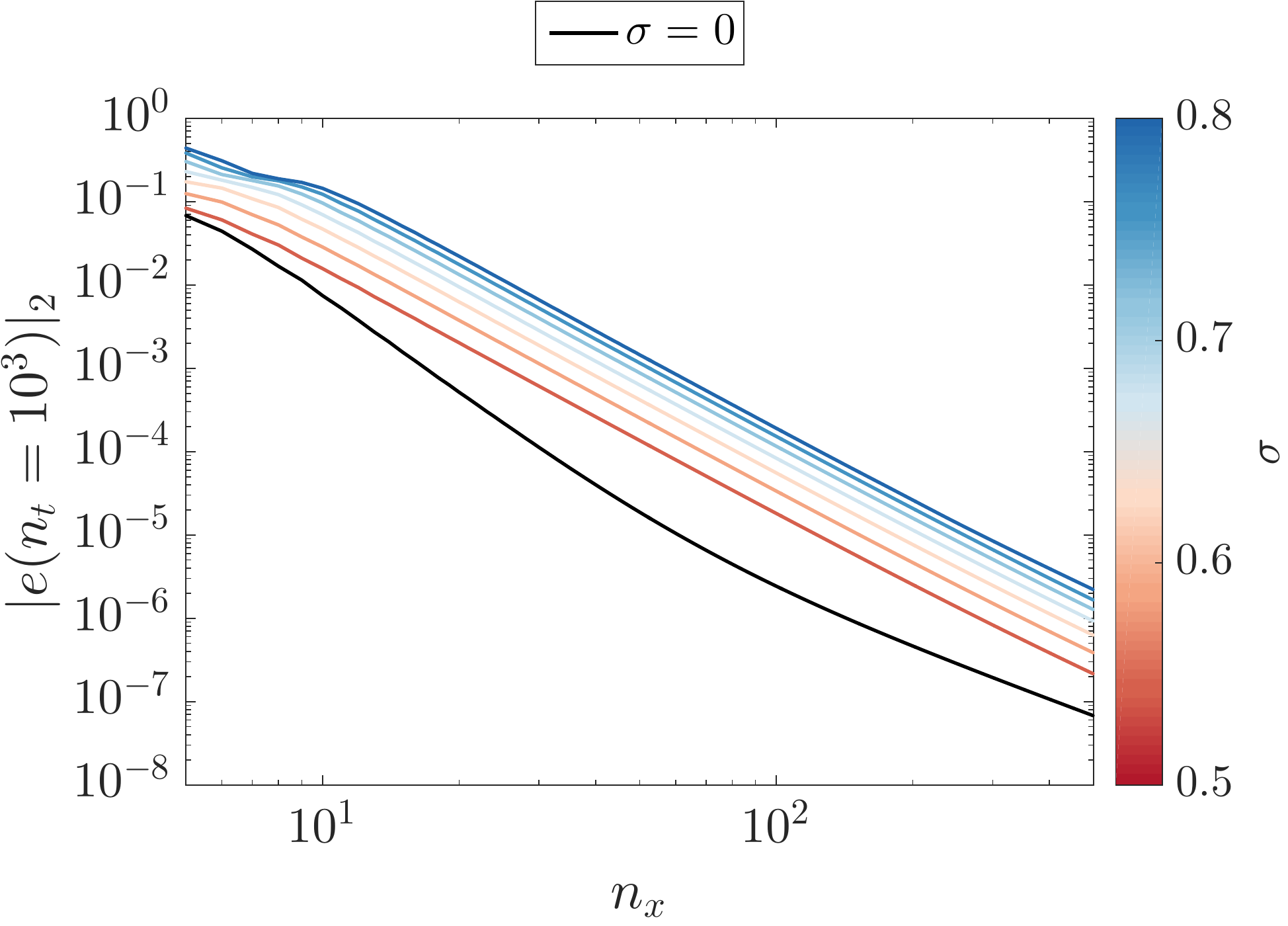}
			\caption{Error calculated of Gaussian bump convection after $10^3$ time iterations.~$p=4$ with Huynh correction, RK33 and constant $\mathrm{CFL}=0.167$, limit for Huynh correction.} 
			\label{fig:FR4_33_bump_error}
		\end{subfigure}		
		~
		\begin{subfigure}[b]{0.45\linewidth}
			\centering
			\includegraphics[width=\linewidth]{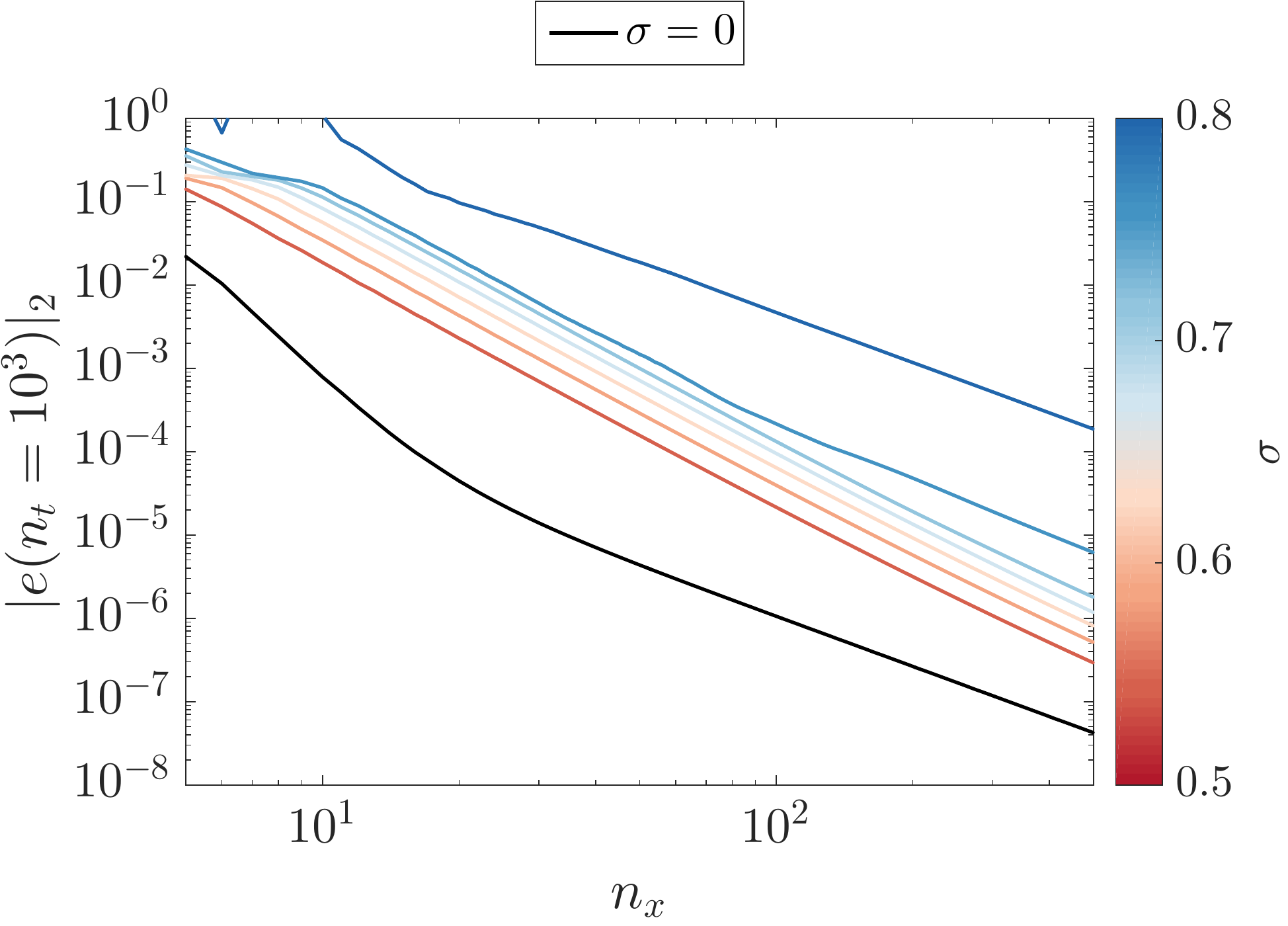}
			\caption{Error calculated of Gaussian bump convection after $10^3$ time iterations.~$p=4$ with Huynh correction, RK44 and constant $\mathrm{CFL}=0.189$, limit for Huynh correction.} 
			\label{fig:FR4_44_bump_error}
		\end{subfigure}
		\caption{Numerical tests of FR using Huynh correction functions, $p=4$.}
		\label{fig:numerical_1}
	\end{figure}

	Validation of the behaviour analytically exhibited is performed via a wave convection test case, in which a wave is input into 1D domain. Interface upwinding is utilised which permits the use of a blow out boundary condition at outlet. During this test three operating conditions are investigated, all of which are at a CFL number that is unstable for the baseline scheme. Two of the waves are for unfiltered schemes and the third has a Gaussian filtered applied with $\sigma=0.6$. The results are displayed in Fig.\ref{fig:FR4_33_k2} and it is clear that application of the filter has stabilised this case, which is in agreement with the earlier findings. 
	
	To highlight the effect of filtering on convective velocity, we present a test in which a Gaussian bump is convected via the linear advection equation through a periodic domain, see Fig.\ref{fig:FR4_33_bump_error}. If the bump is prescribed to be sufficiently narrow this will excite a significant proportion of the wavenumber space, highlighting any instability. The initial condition was therefore chosen to be:
	\begin{equation}
		u = \exp{\big(-(x-0.5)^2/\varsigma^2\big)} 
	\end{equation}    
	where $\varsigma = 0.2$. Comparison of the filtered and unfiltered bump cases clearly shows that the convective velocity has been reduced in the filtered cases, with the diffusive effect of the filter clearly displayed as a reduction in the peak value. This reduction in the convective velocity is consistant with the reduction in the group and phase velocities shown in Fig.\ref{fig:FR4_33_real_170_06}.
	
	 A filter width was also tested that exceeded that the maximum value found to be stable in Section~\ref{sec:vn_results}. This was motivated by the need to understand the mechanism driving the scheme unstable, and it is clear from Fig.\ref{fig:FR4_33_bump} that at $x\approx0.1$ the conserved variable has dipped below zero. Hence, at this level of filtering, the filter is redistributing mass faster than incoming mass can replace it. However, this does seem to be gradient dependant originating from the application of the filter to the differentiation operator, $\mathbf{D}$ and exemplified by the negative derivative on the upwind downwards slope seeming to lead to this region of negative conserved variable.
	
	The analytic solution of this problem permits the calculation of the error which is shown in Fig.~\ref{fig:FR4_33_bump_error}~\&~\ref{fig:FR4_44_bump_error}. The effect of filtering on the error of the scheme has not been significant, causing an upwards shift in the error together with an a small decrease in the order accuracy. This decrease in the order accuracy is to be expected as it has previously been shown by Asthana~\etal~\cite{Asthana2015a} that filtering is equivalent to additionally solving various orders of diffusive derivatives. Although it could be proposed that in order to increase the time step size the grid could be coarsened, in doing so the Nyquist wavenumber will also be reduced and hence the advantage of this method of filtering, due to the already high dissipation at $\hat{k}\approx5\pi/8$, little spectral information is lost and that which is lost, at the very highest wavenumbers, cannot have significant impact on solution due to the impassable dissipation gap seen in Fig.~\ref{fig:FR4_33_t166}~\&~\ref{fig:FR4_33_t170}. What may also be observed are the diverging errors of the unstable filter width, which is consistent with the findings of Fig.\ref{fig:FR4_RK33_1}.

	\subsection{Taylor-Green Vortex}
		Finally, we apply the methodology to an investigation more representative of engineering applications, the Turbulent Taylor-Green Vortex. The primitive variables are initialised as:
	\begin{align}
		u &= U_0\sin{\bigg(\frac{x}{L}\bigg)}\cos{\bigg(\frac{y}{L}\bigg)}\cos{\bigg(\frac{z}{L}\bigg)} \\
		v &= -U_0\cos{\bigg(\frac{x}{L}\bigg)}\sin{\bigg(\frac{y}{L}\bigg)}\cos{\bigg(\frac{z}{L}\bigg)} \\
		w &= 0 \\
		p &= p_0 + \frac{\rho_0U_0^2}{16}\bigg(\cos{\bigg(\frac{2x}{L}\bigg)} + \cos{\bigg(\frac{2y}{L}\bigg)}\bigg)\bigg(\cos{\bigg(\frac{2z}{L}\bigg)} + 2\bigg) \\
		\rho &= \frac{p}{RT_0}
	\end{align}
	The specific case investigated is:
	\begin{equation}		
		R_e = \frac{\rho_0U_0L}{\mu}, \quad\quad
		P_r = 0.71 = \frac{\mu\gamma R}{\kappa(\gamma-1)}, \quad\quad
		M_a = 0.08 = \frac{U_0}{\sqrt{\gamma R T_0}} 
	\end{equation}
	\begin{equation}
		U_0 = 1, \quad\quad \rho_0 = 1, \quad\quad p_0 = 100, \quad\quad R = 1, \quad\quad \gamma = 1.4
	\end{equation}
	This case follows the case set up of Brachet~\etal~\cite{Brachet1983} and DeBonis~\cite{DeBonis2013}. Their DNS data is used as comparison, with the main concern of this investigation being the domain averaged kinetic energy and kinetic energy dissipation.	
	\begin{equation}
		E_k = \frac{1}{2\rho_0\Omega}\int_{\Omega}\rho u_iu_i d\Omega
	\end{equation}		
	
	To isolate the effect of filtering for various various $R_e$, the Taylor-Green test case is used with cell $R_e$ fixed at $R_{e,\mathrm{cell}}=40$. 
	
	\begin{figure}
		\centering
		\begin{subfigure}[b]{0.45\linewidth}
			\centering
			\includegraphics[width=\linewidth]{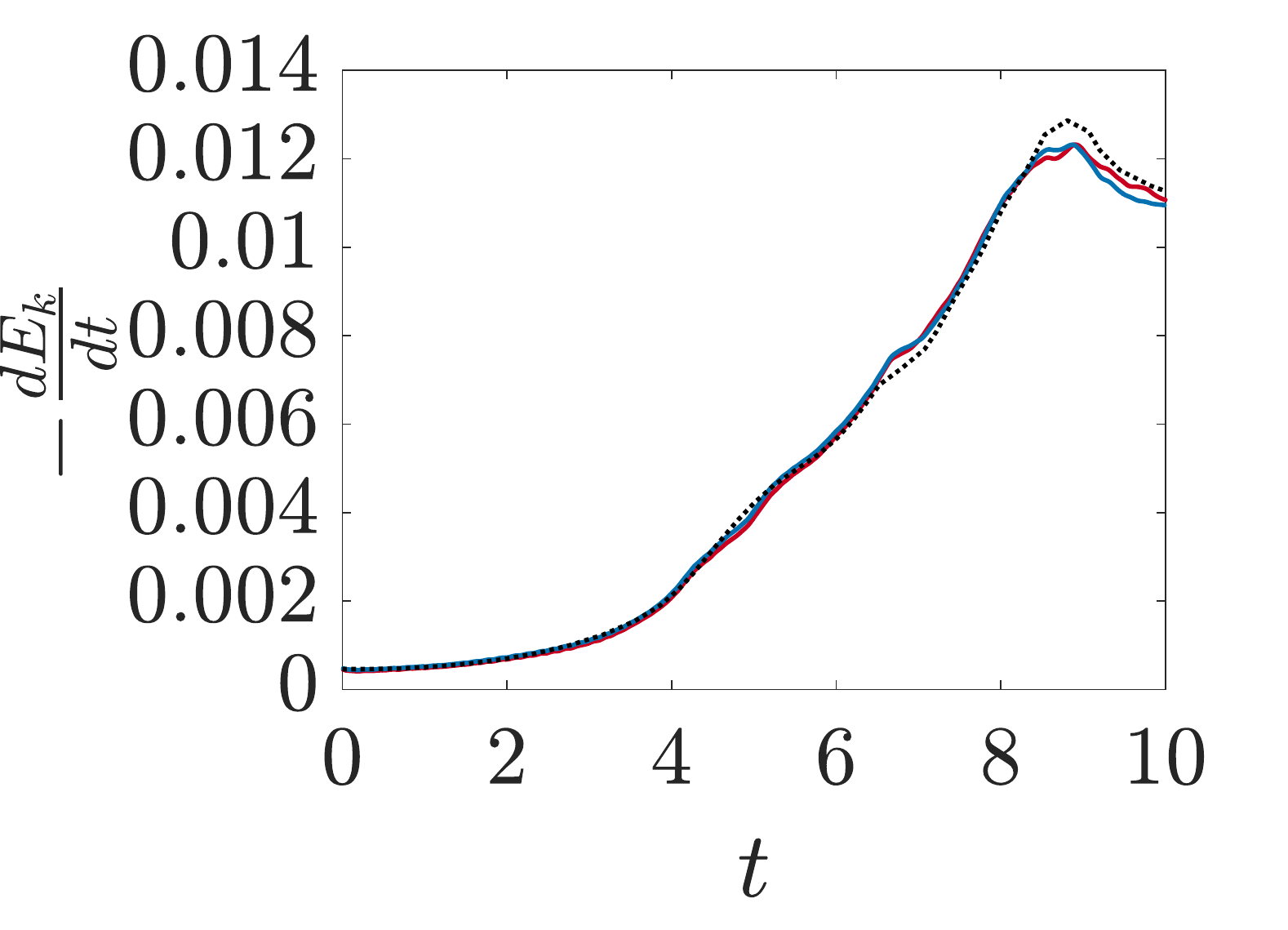}
			\caption{$R_e=1600$, $N=40^3$} 
			\label{fig:FR4_TGV_1600}
		\end{subfigure}
		~
		\begin{subfigure}[b]{0.45\linewidth}
			\centering
			\includegraphics[width=\linewidth]{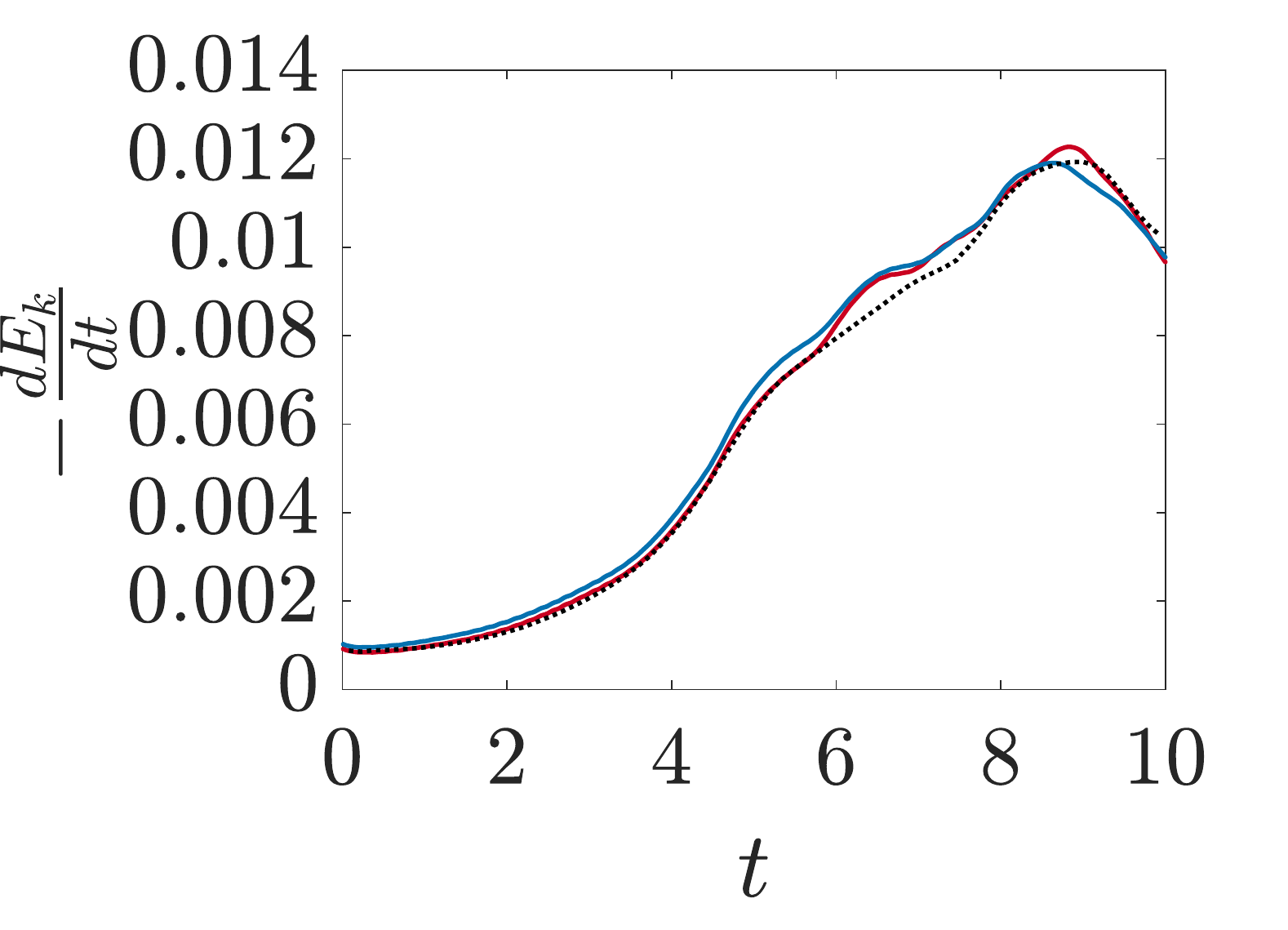}
			\caption{$R_e=800$, $N=20^3$} 
			\label{fig:FR4_TGV_800}
		\end{subfigure}
		~
		\begin{subfigure}[b]{0.45\linewidth}
			\centering
			\includegraphics[width=\linewidth]{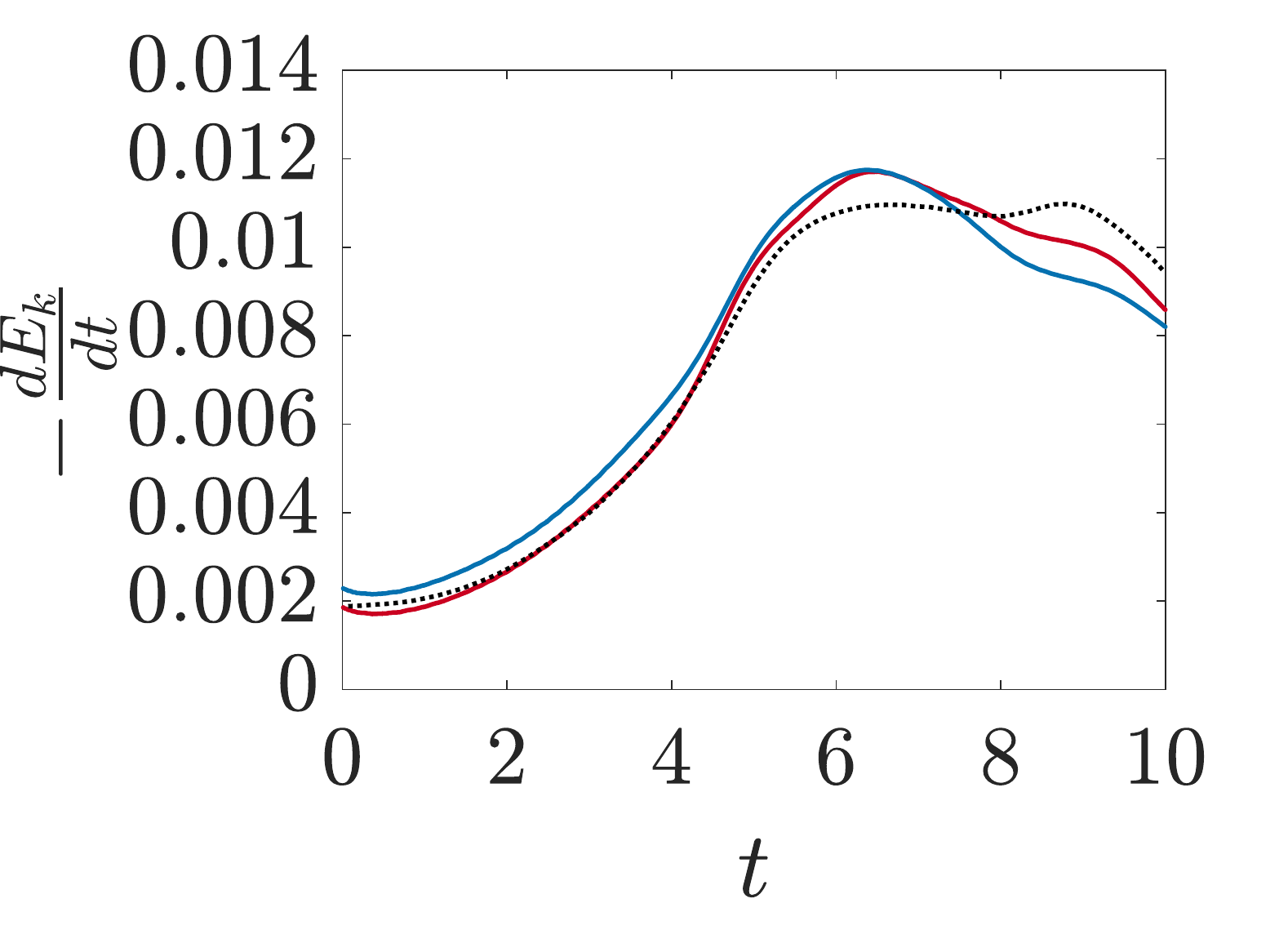}
			\caption{$R_e=400$, $N=10^3$} 
			\label{fig:FR4_TGV_400}
		\end{subfigure}
		~
		\begin{subfigure}[b]{0.45\linewidth}
			\centering
			\includegraphics[width=\linewidth]{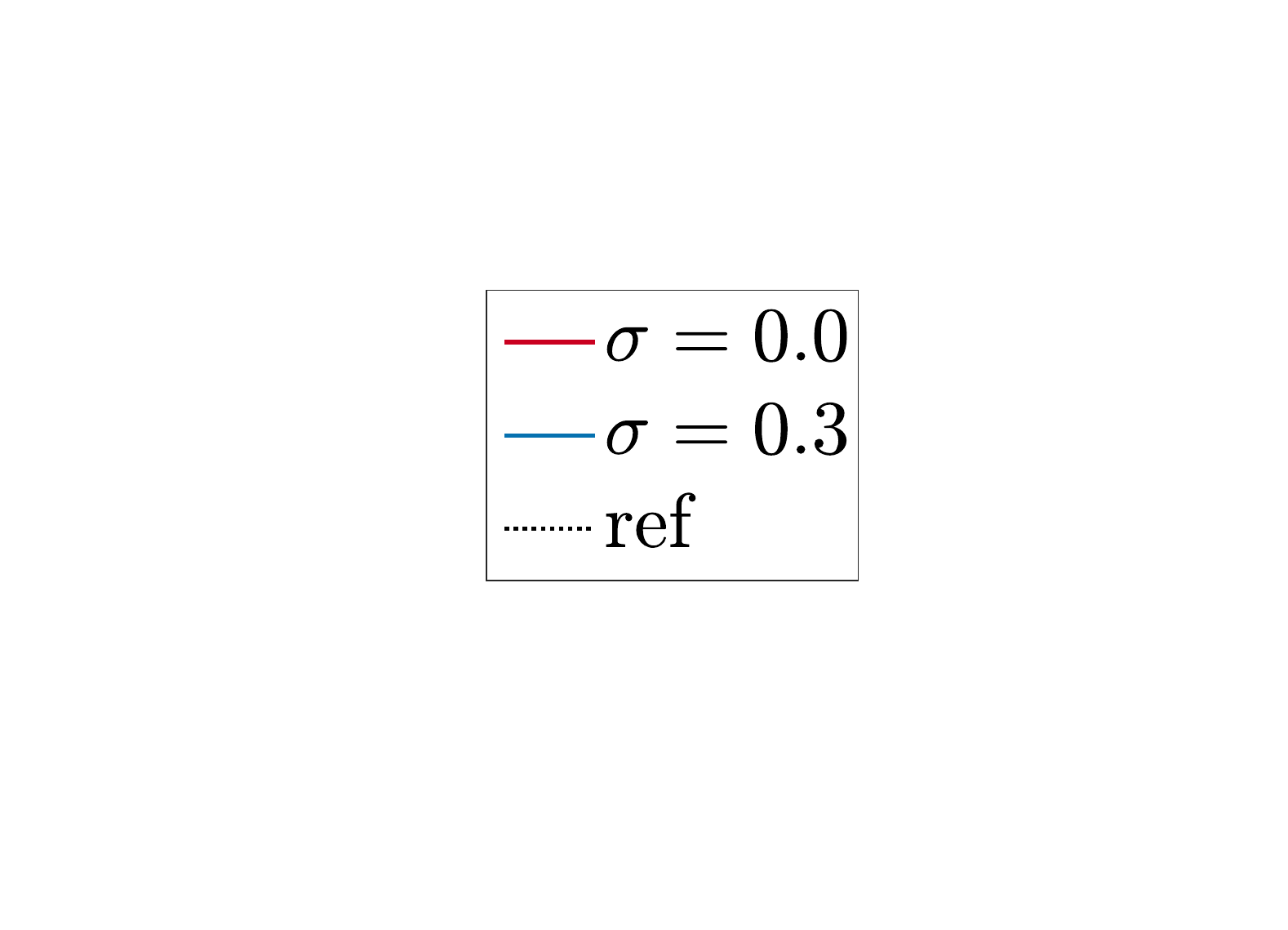}
			\label{fig:FR4_TGV_legend}
		\end{subfigure}
		\caption{Taylor-Green Vortex turbulent kinetic energy dissipation. For $p=4$, using Rusanov flux an the inviscid common interfaces and BR1 at the viscous common interfaces and DG correction functions. The method of temporal integration is RK44 with $dt=1\times10^{-3}$. Between each case the cell $R_e$ was held constant at $R_{e,\mathrm{cell}}=40$. The reference data is that of Brachet~\etal~\cite{Brachet1983}}
		\label{fig:FR4_TGV_40}
	\end{figure}
	
	What should be clear from Fig.\ref{fig:FR4_TGV_1600} is that at this resolution the transition to turbulence is well captured , with only minor deviation from peak kinetic energy dissipation. This is a result that is shared by  Bull and Jameson~\cite{Bull2014a}. However, as we move to lower $R_e$, Fig.\ref{fig:FR4_TGV_400}, the deviations grow, which may be explained by the increased sensitivity at lower $R_e$ to changes in the apparent viscosity. With a filter added the effect on lower $R_e$ cases also becomes more pronounced, this is again due to the increased sensitivity of lower $R_e$ flows to changes in apparent viscosity. However, there is another effect that impacts the results here, that the the filter spectral width is constant in the reference domain and so depends on the physical size of the element. The result of this is that as the grid becomes coarser the filter becomes more dissipative, \emph{i.e.} the filter $R_e$ depends inversely on element length. From the filter viscosity of derived by Asthana~\etal~\cite{Asthana2015a} and the filter definition of Eq.(\ref{eq:gaussian_filter}) the filter $R_e$ can be defined as:
	\begin{equation}
		R_{e,\mathrm{filter}} = \frac{24\rho u\tau}{\sigma^2h}	
	\end{equation}
	where $\tau$ is the time step size and $h$ is the grid spacing. To investigate further the effect of grid spacing and filter width on the simulation of the full Navier-Stokes equations, a series of different grid spacings and filter width are considered for constant global $R_e=1600$.
	
	\begin{figure}
		\centering
		\begin{subfigure}[b]{0.7\linewidth}
			\centering
			\includegraphics[width=\linewidth]{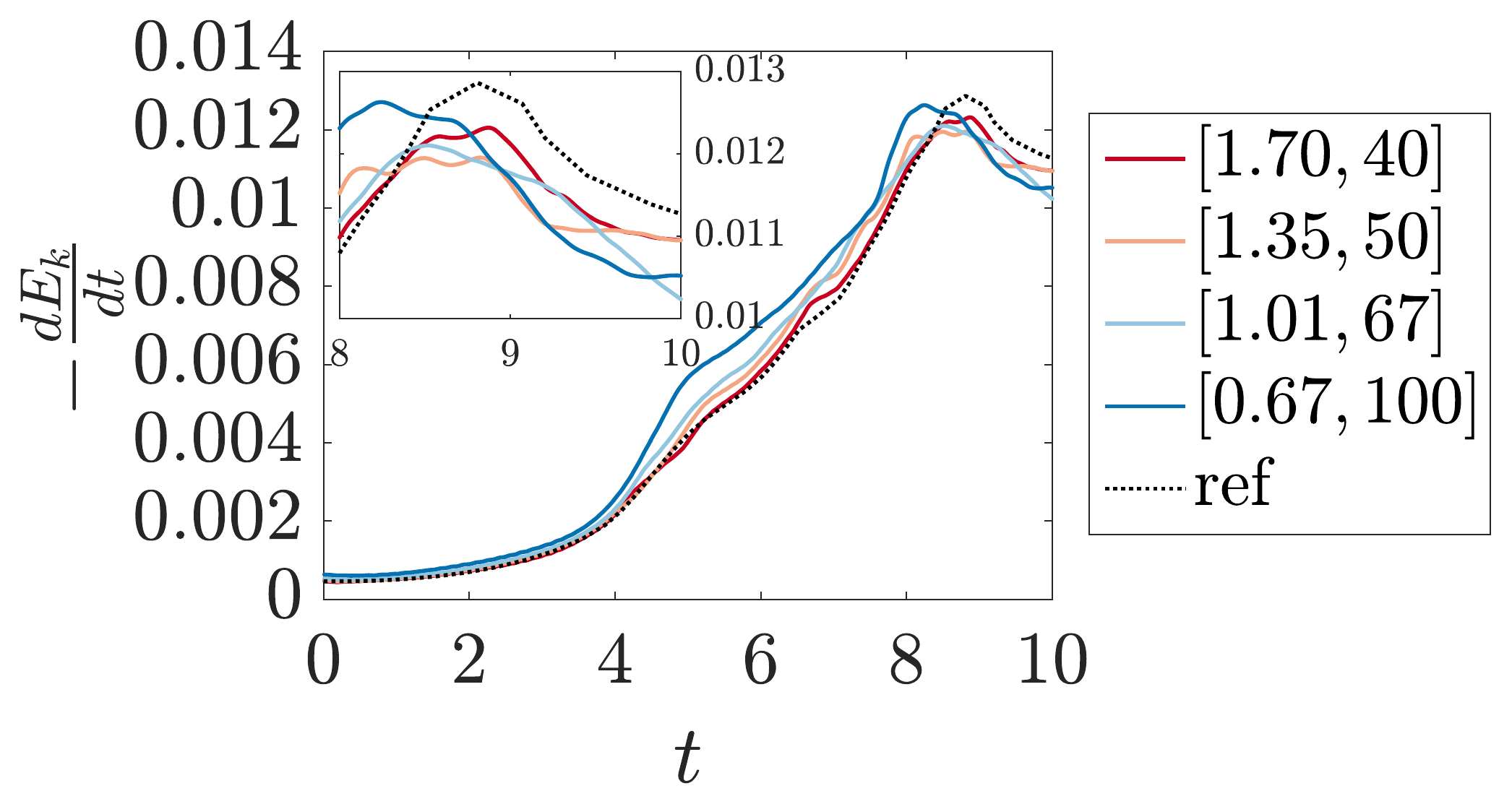}
			\caption{$R_e = 1600$ on grids $40^3$, $32^3$, $24^3$, and $16^3$ with $\sigma=0.3$.} 
			\label{fig:FR4_TGV_1600_dx}
		\end{subfigure}
		~
		\begin{subfigure}[b]{0.7\linewidth}
			\centering
			\includegraphics[width=\linewidth]{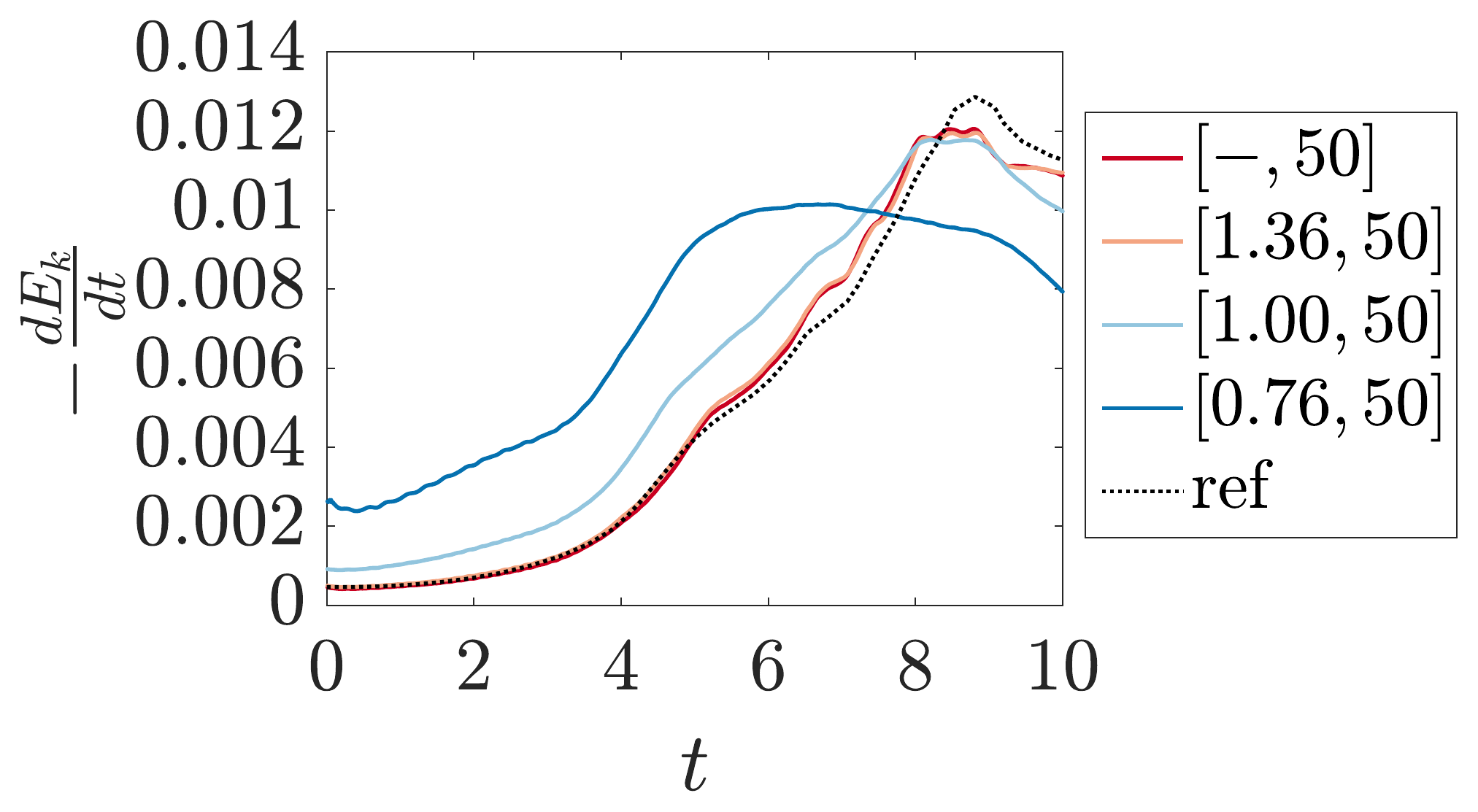}
			\caption{$R_e = 1600$ with $32^3$ elements with $\sigma\in\{0,0.3,0.35,0.4\}$.} 
			\label{fig:FR4_TGV_1600_ds}
		\end{subfigure}
		\caption{Turbulent Taylor-Green vortex kinetic energy dissipation for FR, $p=4$, using nodal DG corrections, with Rusanov inviscid and BR1 viscous common interface flux calculations. The temporal integration method is RK44 with $\tau=1\times10^{-3}$. The legend format is $[R_{e,\mathrm{filter}},R_{e,\mathrm{cell}}]$. The reference data is that of Brachet~\etal~\cite{Brachet1983}}
		\label{fig:FR4_TGV_1600_dxs}
	\end{figure}
	
	What can be seen from Fig.\ref{fig:FR4_TGV_1600_dxs} is that there is a marked impact of the filtering as $R_{e,\mathrm{filter}}$ reduces below one, with exact transition dependent on the resolution of the grid relative to the flow. Furthermore, the definition of filter viscosity used in the definition of $R_{e,\mathrm{filter}}$ takes into account only the lowest order diffusion derivative and hence there will be some higher order dependence of filter Reynolds number on grid spacing that will also affect the exact point of transition. The physical interpretation of this result is that as the filter $R_e$ decreases, the rate at which solutions are convected through an element reduces compared to the time scale of the diffusion caused by the filter. Therefore, a reduction of $R_{e,\mathrm{filter}}$ below one implies that the time scale of diffusion has become more significant than that of inertia, and so the filter diffusion becomes far more apparent. 
	
\section{Conclusions}
	It has been shown that filtering in the reference domain can be used to alter the behaviour of the underlying FR scheme by the addition of dissipation at higher wavenumbers. This dissipation can be used to increase the temporal stability, which, for the case of $p=4$ with RK33 temporal integration, can increase the CFL limit by $\approx25\%$. For temporal stabilisation the most successful method of filtering FR is complete filtering of the spatial scheme, as it provides a more coherent treatment of higher wavenumbers. Hence, this form of filtering in the reference domain can be thought of as modification of the numerics for implicit LES, rather than filtering for explicit LES.
	
	Numerical tests were used to validate analytical results, with an investigation of order accuracy showing that filtering only mildly reduces the order of accuracy of the method. Finally, a filter Reynolds number is defined and, through investigation with turbulent flows, filtering is shown to have a significant impact when the filter Reynolds number reduces below a critical value, which in this case that was found to be one.

\section*{Acknowledgments}
The support of the Engineering and Physical Sciences Research Council of the United Kingdom is gratefully acknowledged.

\bibliographystyle{aiaa}
\bibliography{library}

\end{document}